\documentclass[reqno,11pt]{amsart}

\usepackage{a4wide}
\usepackage[top=0.5in,bottom=0.7in]{geometry}
\usepackage{graphicx}
\counterwithin{figure}{section} 
\usepackage{bm}
\usepackage{xcolor}
\usepackage{imakeidx}
\usepackage{tikz}
\usepackage{caption}
\usepackage[bookmarks=false]{hyperref}
\usepackage{amsmath,amsthm,amssymb}
\usepackage{mathtools}
\usepackage[justification=raggedright]{caption}
\usepackage{xspace}
\usepackage{complexity}
\usepackage{bbold}
\usepackage{enumitem}
\usepackage{physics}
\usepackage{derivative}
\usepackage{mathrsfs}
\usepackage{cleveref}
\usepackage{float}
\usepackage{relsize}

\makeindex[title={Index\label{index}},columns=2]

\newcommand{\mathdf}[1]{#1\indexmath{#1}}
\def\indexmath#1{\index{\ensuremath{#1}}}

\numberwithin{equation}{section}
\newtheorem{theorem}{Theorem}[section]
\newtheorem{lemma}[theorem]{Lemma}
\newtheorem{corollary}[theorem]{Corollary}
\newtheorem{assumption}[theorem]{Assumption}
\newtheorem{proposition}[theorem]{Proposition}
\theoremstyle{definition}
\newtheorem{definition}[theorem]{Definition}
\theoremstyle{remark}
\newtheorem{remark}[theorem]{Remark}

\def\cW{\mathcal{W}}
\def\hW{\widehat{\mathcal{W}}}
\def\cH{\mathcal{H}}
\def\cS{\mathcal{S}}
\def\cN{\mathcal{N}}
\newcommand\cT{{\mathcal{T}}}
\newcommand\hT{{\widehat{\mathcal{T}}}}
\DeclareMathOperator{\Int}{Int}
\newcommand{\half}{\frac{1}{2}}

\newcommand{\online}{}  

\title[The discrete Laplace method and the 3XOR-SAT problem]{The discrete Laplace asymptotic method and its application to the 3XOR satisfiability problem}

\author{J.~William~Helton}
\address{(JWH) University of California, Department of Mathematics, San Diego, USA}
\email{whelton@ucsd.edu}

\author{Jared~A.~Hughes}
\address{(JAH) University of California, Department of Mathematics, San Diego, USA}
\email{jahughes241@gmail.com}

\author{Peter~Schlosser}
\address{(PS) Graz University of Technology, INstitute of Applied Mathematics, Austria}
\email{pschlosser@math.tugraz.at}

\begin{document}

\begin{abstract}
A standard way to calculate the asymptotic behavior of integrals of the form $\int_\mathcal{W}g(x)e^{-nh(x)}dx$ is the (continuous) Laplace asymptotic method. However, also discrete sums like $\sum_{x\in\mathcal{W}\cap\Lambda_n}g_n(x)e^{-nh_n(x)}$ have similar behavior, when $\Lambda_n$ is a discrete grid which becomes infinitely fine, and the functions $g_n$ and $h_n$ converge to $g$ and $h$ respectively. We go even further, and also derive the asymptotic formula for sums of the form $\sum_{x\in\mathcal{W}\cap\Lambda_n}\mathcal{S}_n(x)$, where the summand $\mathcal{S}_n (x)$ asymptotically behaves as $g_n(x)e^{-nh_n(x)}$. 

The motivation, and also an immediate application, will be filling in all details in the classical breakthrough paper \cite{dBM-FOCS}, which gives the solvability (phase transition) threshold of the 3XOR-SAT problem 
using the second moment method. 
Various analytical arguments there were lightly described, but the \online appendix to this paper combines recent results to fill all of them in. We would expect our theorems on asymptotics to apply to other (especially combinatorial) problems as well. For example, they seem effective on 3XOR-GAME problems.
\end{abstract}


\maketitle

\textbf{Acknowledgements}: Peter Schlosser was funded by the Austrian Science Fund (FWF) under Grant No. J 4685-N and by the European Union--NextGenerationEU.

\section{Introduction}

The basic problem the discrete Laplace asymptotic method
addresses goes as follows.
The setting is in 
 a bounded region $\cW \subseteq\mathbb{R}^d$, for a fixed 
integer $d\geq 1$.
Moreover, we consider a sequence of lattices
\index{lattice, $\Lambda_k$}
\index{$\Lambda_k, A_k,  n_k,  v_k$}
\index{$A,v$}
\begin{equation}\label{Eq_Lattice}
\Lambda_k\coloneqq\frac{1}{n_k}A_k\mathbb{Z}^d+v_k,
\end{equation}
where $(n_k)_k\in\mathbb{N}$ is a sequence of integers with $n_k\geq k$, $(v_k)_k\in\mathbb{R}^d$ a sequence of vectors with $\lim_{k\to\infty}v_k=v\in\mathbb{R}^d$, and $(A_k)_k\in\mathbb{R}^{d\times d}$ is a sequence of matrices which converge as $\lim_{k\to\infty}A_k=A$ to an invertible matrix $A\in\mathbb{R}^{d\times d}$. Given 
$\cS_k:\cW\cap\Lambda_k\to\mathbb{C}$ a sequence of functions with asymptotic approximation $\cS_k(x)\sim g_k(x)e^{-n_kh_k(x)}$, we want to see how 
\begin{align*}
    \sum_{x\in\cW\cap\Lambda_k}\cS_k(x)
\end{align*}
behaves asymptotically as $k\to\infty$. If we are in a case where $g_k\to g$ and $h_k\to h$, the sum may behave like a Riemann sum, and we optimistically obtain
\begin{align*}
\sum_{x\in\cW\cap\Lambda_k}\cS_k(x)\sim\sum_{x\in\cW\cap\Lambda_k}g_k(x)e^{-n_kh_k(x)}\sim\frac{n_k^d}{|\det(A_k)|}\int_\cW g_ke^{-n_kh_k}dx\sim\Omega_k,
\end{align*}
where the classical Laplace method for integrals applies to give
\begin{align*}
    \Omega_k\coloneqq\frac{(2\pi n_k)^{\frac{d}{2}}g(x_0)e^{-n_kh(x_0)}}{|\det(A)|\sqrt{\det(\mathcal{H}\{h\}(x_0))}}.
\end{align*}
Here, $x_0\in\cW$ is the point where the global minimum of the function $h$ is obtained, and $\mathdf{\mathcal{H}\{h\}(x_0)}$ is the Hessian matrix evaluated at this point. 
The goal of this paper is to give reasonably 
general conditions under which these results do hold. More precisely, our first main result, \Cref{thm_Discrete_Laplace_method}, implies that
\begin{align*}
\sum_{x\in\cW\cap\Lambda_k}g_k(x)e^{-n_kh_k(x)}\sim\Omega_k
\end{align*}
holds under precise conditions. Our second main result, \Cref{thm_Discrete_Sk}, then even proves
\begin{equation*}
\sum_{x\in\cW\cap\Lambda_k}\cS_k(x)\sim\Omega_k,
\end{equation*}
in the same setting, plus a few assumptions connecting $\cS_k$ to $g_ke^{-n_kh_k}$. \medskip

The papers' final section, \Cref{sec_Application_3XORSAT_3XORGAME}, describes two applications to special systems of linear equations over the binaries; these correspond to classical 3XOR satisfiability problems (3XOR-SAT) and to 3XOR games (3XOR-GAME). The issue is the ratio of number of equations $m$ to number of unknowns $n$ for which the system of equations is likely to be solvable: Does there exist a ``solvability threshold" $c^*$, so that asymptotically $m/n<c^*$ implies solvable, and $m/n>c^*$ implies not solvable, with high probability?

The classic beautiful tour de force 
of Dubois--Mandler \cite{dBM-FOCS} says yes for 3XOR-SAT and computes $c^*$. 
One of the many steps in their proof is \cite[Lemma 3.6]{dBM-FOCS}, which asserts a discrete Laplace asymptotic formula,
and gives the intuitive justification that it is analogous to the classical integral Laplace asymptotic formula.

\Cref{sec_Application_3XORSAT_3XORGAME} rigorously proves their conclusion as an application of our theorems.
In support of this, and also using recent results of \cite{stirlingConnDobro,CHHarxiv}, we produced an \online appendix, which fills in all of the analytical details which were alluded to, but not included in \cite{dBM-FOCS}.

After this, in \Cref{sec:3XOR-GAME},  we point out that analysis of the 3XOR-GAME thresholds 
produce a big messy sum whose asymptotics our theory addresses. In unpublished work we have filled in most of the details.

\section{The discrete Laplace asymptotic method}

This section gives and proves a discrete variation of the Laplace asymptotic method. As a motivation, in \Cref{lemma:laplace-2} it starts with the (already known) continuous analog of what we are proving. Discretizing the domain of integration then turns the continuous Laplace method into the discrete Laplace method, see \Cref{thm_Discrete_Laplace_method}. Finally, in \Cref{thm_Discrete_Sk}, we further generalize our results by replacing the function $g_ke^{-n_kh_k}$ by arbitrary discrete summands $S_k$, which only asymptotically behave as $g_ke^{-n_kh_k}$. \medskip

\subsection{The continuous Laplace method}

To give context for the situation, we start with the classical continuous Laplace method, which our discrete Laplace method in \Cref{sec_Discrete_Laplace_method} will imitate, and which the asymptotics for discrete summands in \Cref{sec_Discrete_summands} extends. The following \Cref{lemma:laplace-2} is a slight variation of the statement in \cite[Theorem~15.2.2]{SimonComplexAnalysis}. For further terms in the respective asymptotic expansion, see \cite[Theorem~15.2.5]{SimonComplexAnalysis}, which relies on $g,h$ being infinitely often differentiable.

\begin{theorem}[Continuous Laplace method]\label{lemma:laplace-2}
Let $\cW\subseteq\mathbb{R}^d$ be open and bounded, and $h\colon\cW\to\mathbb{R}$. Let us assume that there exists some $x_0\in\cW$, such that $h\in C^2(\cN_0)$ in an open neighborhood $\cN_0$ of $x_0$, and there holds the following two assumptions: \medskip

\begin{enumerate}
\item[i)] The Hessian $\cH\{h\}(x_0)$ is strictly positive definite. \medskip

\item[ii)] For every open neighborhood $\cN$ of $x_0$, there is $\inf_{x\in\cW\setminus\cN}h(x)>h(x_0)$. \medskip
\end{enumerate}

Then for every Riemann integrable function $g\colon\cW\to\mathbb{C}$, which is continuous near $x_0$, there converges
\begin{equation}\label{Eq_Continuous_Laplace_method}
\lim\limits_{n\to\infty}n^{\frac{d}{2}}e^{nh(x_0)}\int_\cW g(x)e^{-nh(x)}dx=\frac{(2\pi)^{\frac{d}{2}}g(x_0)}{\sqrt{\det(\cH\{h\}(x_0))}}.
\end{equation}
\end{theorem}

\begin{remark}
Before we start the Discrete Laplace Method, we emphasize that integrability over all of $\cW$ is not significant, since the Discrete Laplace Method concerns sums over a lattice intersected with $\cW$. Instead, we only need integrability over a small neighborhood $\cN_0$ where $h$ is minimized, and relevant bounds outside the neighborhood.
\end{remark}

\subsection{The Discrete Laplace method}\label{sec_Discrete_Laplace_method}

The idea of the discrete Laplace method is to approximate the integral on the left hand side of \eqref{Eq_Continuous_Laplace_method} by its Riemann sum. In other words, we will partition the integration domain $\cW$ into a lattice $\Lambda_k$ of the form \eqref{Eq_Lattice}, and approximate
\begin{equation}\label{Eq_Discrete_Laplace_formal_asymptotics}
\int_\cW g(x)e^{-n_k h(x)}dx\sim\frac{\abs{\det(A_k)}}{n_k^d}\sum\limits_{x\in\cW\cap\Lambda_k}g_k(x)e^{-n_kh_k(x)},
\end{equation}
where the prefactor $\frac{|\det(A_k)|}{n_k^d}$ is the cell volume of the lattice. In \Cref{thm_Discrete_Laplace_method} we then show that under the following \Cref{ass_Laplace_assumption_discrete}, the right hand side of this approximation still converges to the same limit as the continuous Laplace method in Eq.\eqref{Eq_Continuous_Laplace_method} suggests.

\begin{assumption}\label{ass_Laplace_assumption_discrete}
Let $\mathdf{\cW}\subseteq\mathbb{R}^d$ be bounded, and for every $k\in\mathbb{N}$ consider a lattice $\Lambda_k$ as in Eq.\eqref{Eq_Lattice}, as well as functions $g_k\colon\cW\cap\Lambda_k\to\mathbb{C}$ and $h_k\colon\cW\to\mathbb{R}$. 
\index{$g_k(x)$}
\index{$h_k(x)$}
Then we say that these objects satisfy the \textnormal{discrete Laplace assumptions}, if there exists some $x_0\in\Int(\cW)$, such that $h_k\in C^2(\cN_0)$ in an open neighborhood $\mathdf{\cN_0}$ of $\mathdf{x_0}$, and there holds the following four assumptions: \medskip

\begin{enumerate}
\item[i)] There exist functions $g\in C(\cN_0)$ and $h\in C^2(\cN_0)$, such that there converges
\index{$g(x)$}
\index{$h(x)$}
\begin{equation*}
\lim\limits_{k\to\infty}\Vert g_k-g\Vert_{\infty,\,\cN_0\cap \Lambda_k}=0\qquad\text{and}\qquad\lim\limits_{k\to\infty}\Vert h_k-h\Vert_{C^2(\cN_0)}=0;
\end{equation*}

\item[ii)] In the point $x_0$ there is
\begin{equation*}
h_k(x_0)=h(x_0),\quad\nabla h_k(x_0)=0,\quad\text{and}\quad\mathcal{H}\{h\}(x_0)\text{ is strictly positive definite};
\end{equation*}

\item[iii)] For every open neighborhood $\cN$ of $x_0$, there is
\begin{equation}\label{Eq_hk_lower_bound}
\liminf_{k\to\infty}\inf_{x\in\cW\setminus\cN}h_k(x)>h(x_0);
\end{equation}

\item[iv)] There exists some $C_g\geq 0$ and $l\in\mathbb{N}$, such that
\begin{equation}\label{Eq_gk_upper_bound_poly}
\Vert g_k\Vert_{\infty,\mathcal{W}\cap\Lambda_k}\leq\mathdf{C_g}k^{\ell},\qquad k\in\mathbb{N}.
\end{equation}
\end{enumerate}
\end{assumption}

The following three Lemmas \ref{lem_Sum_outside}, \ref{lem_Sum_Gaussian} and \ref{lem_Taylor_estimate} will now be preparatory for the main theorem of this subsection, Theorem \ref{thm_Discrete_Laplace_method}. In \Cref{lem_Sum_outside} we will prove that the sum in Eq.\eqref{Eq_Discrete_Laplace_formal_asymptotics} tends to zero on every region $\cW\setminus\cN$ away from $x_0$, where $h$ is strictly larger than $h(x_0)$, see Eq.\eqref{Eq_hk_lower_bound}.

\begin{lemma}\label{lem_Sum_outside}
Let $\cW\subseteq\mathbb{R}^d$ be bounded, and for every $k\in\mathbb{N}$ consider a lattice $\Lambda_k$ as in Eq.\eqref{Eq_Lattice}, as well as functions $g_k\colon\mathcal{W}\cap\Lambda_k\to\mathbb{C}$ and $h_k\colon\mathcal{W}\to\mathbb{R}$, which satisfy \Cref{ass_Laplace_assumption_discrete}. Then for every open neighborhood $\mathdf{\cN}\subseteq \cW$ of $x_0$, there converges
\begin{equation*}
\lim\limits_{k\to\infty}n_k^{-\frac{d}{2}}e^{n_kh(x_0)}\sum\limits_{x\in(\cW\setminus\cN)\cap\Lambda_k}g_k(x)e^{-n_kh_k(x)}=0.
\end{equation*}
\end{lemma}

\begin{proof}
Since $\cW$ is bounded, and the cell volume of the lattice in Eq.\eqref{Eq_Lattice} is $\abs{\det(A_k)}n_k^{-d}$, with $\abs{\det(A_k)}\to\abs{\det(A)}$, 
there exists a constant $C_A\geq 0$, such that the number of lattice points inside $\cW$ is bounded by
\begin{equation*}
\abs{\cW\cap\Lambda_k}\leq\mathdf{C_A}n_k^d,\qquad k\in\mathbb{N}.
\end{equation*}
Since the limit inferior in Eq.\eqref{Eq_hk_lower_bound} is strictly larger than $h(x_0)$, there exists some $N\in\mathbb{N}$, such that
\begin{equation*}
\mathdf{h_{\min}}\coloneqq\inf_{k\geq N}\inf_{x\in\cW\setminus\cN}h_k(x)>h(x_0).
\end{equation*}
Using also the upper bound \eqref{Eq_gk_upper_bound_poly} of the norm $\Vert g_k\Vert_{\infty,\mathcal{W}\cap\Lambda_k}$, as well as $n_k \geq k$, there converges
\begin{align*}
n_k^{-\frac{d}{2}}e^{n_kh(x_0)}\sum\limits_{x\in(\cW\setminus\cN)\cap\Lambda_k}\!\abs{g_k(x)}\,e^{-n_kh_k(x)}&\leq n_k^{-\frac{d}{2}}\big|(\mathcal{W}\setminus\mathcal{N})\cap\Lambda_k)\big|C_gk^le^{-n_k(h_{\min}-h(x_0))} \\
&\leq C_A \,C_g \,n_k^{\frac{d}{2} }n_k^\ell\,e^{-n_k(h_{\min}-h(x_0))} \overset{k\rightarrow\infty}{\longrightarrow}0. \qedhere
\end{align*}
\end{proof}

The next lemma now handles the sum \eqref{Eq_Discrete_Laplace_formal_asymptotics} in neighborhoods $\cN$ of $x_0$, and with a function $h_k$ in the exponent which is of the form $\frac{1}{2}(x-x_0)^\top H(x-x_0)$. \Cref{thm_Discrete_Laplace_method} will apply this 
lemma to a small ball around $x_0$, where $h$ is locally quadratic.

\begin{lemma}\label{lem_Sum_Gaussian}
Let $H\in\mathbb{R}^{d\times d}$ be strictly positive definite, $x_0\in\mathbb{R}^d$, and $(\Lambda_k)_k$ a sequence of lattices of the form \eqref{Eq_Lattice}. Then for every open neighborhood $\cN$ of $x_0$, there converges
\begin{equation}\label{Eq_Sum_Gaussian}
\lim\limits_{k\to\infty}n_k^{-\frac{d}{2}}\sum_{x\in\cN\cap\Lambda_k}e^{-\frac{n_k}{2}(x-x_0)^\top H(x-x_0)}=\frac{(2\pi)^{\frac{d}{2}}}{|\det(A)|\sqrt{\det(H)}}.
\end{equation}
\end{lemma}

\begin{proof}
Applying the change of variables $y:=\sqrt{n_k}\,(x-x_0)$, transforms the lattice $\Lambda_k$ and the neighborhood $\mathcal{N}$ into
\begin{equation}\label{Eq_Lattice_tilde}
\mathdf{\widetilde{\Lambda}_k}:=\frac{1}{\sqrt{n_k}}A_k\mathbb{Z}^d+\sqrt{n_k}\,(v_k-x_0),\qquad\text{and}\qquad\mathdf{\mathcal{N}_k}:=\sqrt{n_k}\,(\mathcal{N}-x_0).
\end{equation}
Under this change of variables, we can rewrite the sum on the left hand side of Eq.\eqref{Eq_Sum_Gaussian} as
\begin{equation}\label{Eq_Sum_Gaussian_1}
\sum_{x\in\cN\cap\Lambda_k}e^{-\frac{n_k}{2}(x-x_0)^\top H(x-x_0)}=\sum\limits_{y\in\widetilde{\cN}_k\cap\widetilde{\Lambda}_k}e^{-\frac{1}{2}y^\top Hy}.
\end{equation}
Next we note that the lattice $\widetilde{\Lambda}_k$ in Eq.\eqref{Eq_Lattice_tilde} has a unit cell volume of $|\det(A_k)|n_k^{-\frac{d}{2}}$, and the neighborhoods $\widetilde{\cN}_k$ cover the whole space $\mathbb{R}^d$ when $k\to\infty$. This means, if we multiply the equation \eqref{Eq_Sum_Gaussian_1} with $|\det(A_k)|n_k^{-d/2}$, it becomes an approximating sum of the Gaussian integral
\begin{align*}
\lim\limits_{k\to\infty}|\det(A_k)|n_k^{-\frac{d}{2}}\sum\limits_{x\in\cN\cap\Lambda_k}e^{-\frac{n_k}{2}(x-x_0)^\top H(x-x_0)}&=\lim\limits_{k\to\infty}|\det(A_k)|n_k^{-\frac{d}{2}}\sum\limits_{y\in\widetilde{\cN}_k\cap\widetilde{\Lambda}_k}e^{-\frac{1}{2}y^\top Hy^\top} \\
&=\int_{\mathbb{R}^d}e^{-\frac{1}{2}y^\top Hy}dy=\frac{(2\pi)^{d/2}}{\sqrt{\det(H)}}.
\end{align*}
Using also $\lim_{k\to\infty}\abs{\det(A_k)}=\abs{\det(A)}\neq 0$, gives the stated limit \eqref{Eq_Sum_Gaussian}.
\end{proof}

The final preparatory \Cref{lem_Taylor_estimate}, gives a quadratic Taylor estimate from above and from below of the functions $h_k$ in a neighborhood of its minimum $x_0$.

\begin{lemma}\label{lem_Taylor_estimate}
In the setting of Assumption~\ref{ass_Laplace_assumption_discrete}, for every $\varepsilon>0$ there exists an open neighborhood $\mathdf{\cN_\varepsilon}$ of $x_0$ and some $\mathdf{K_\varepsilon}\in\mathbb{N}$, such that for every $x\in\cN_\varepsilon$, $k\geq K_\varepsilon$
\begin{equation}\label{Eq_Taylor_estimate}
\frac{1-\varepsilon}{2}(x-x_0)^\top\cH\{h\}(x_0)(x-x_0)\leq h_k(x)-h(x_0)\leq\frac{1+\varepsilon}{2}(x-x_0)^\top\cH\{h\}(x_0)(x-x_0).
\end{equation}
\end{lemma}

\begin{proof}
Without loss of generality we will only consider $x_0=0$, $h(x_0)=0$, and assume that the neighborhood $\cN_0$ in \Cref{ass_Laplace_assumption_discrete} is convex. Since $h_k\in C^2(\cN_0)$ is twice continuously differentiable, with $h_k(0)=0$ and $\nabla h_k(0)=0$ by Assumption~\ref{ass_Laplace_assumption_discrete}~ii), the Taylor formula provides for every $x\in\cN_0$ some $\xi_{k,x}\in[0,x]$, such that
\begin{equation}\label{Eq_Taylor_estimate_1}
h_k(x)=\frac{1}{2}x^\top\cH\{h_k\}(\xi_{k,x})x,\qquad x\in\cN_0.
\end{equation}
Since there is also $\cH\{h\}(0)$ strictly positive definite, again by Assumption~\ref{ass_Laplace_assumption_discrete}~ii), there exists some $h_0>0$, such that
\begin{equation*}
x^\top\cH\{h\}(0)x\geq h_0|x|^2,\qquad x\in\mathbb{R}^d\setminus\{0\}.
\end{equation*}
Consequently, we can estimate
\begin{equation*}
\frac{\big|x^\top\big(\cH\{h_k\}(\xi_{k,x})-\cH\{h\}(0)\big)x\big|}{x^\top\cH\{h\}(0)x}\leq\frac{\big\Vert\cH\{h_k\}(\xi_{k,x})-\cH\{h\}(0)\big\Vert}{h_0},\qquad x\in\mathbb{R}^d\setminus\{0\}.
\end{equation*}
This inequality now turns the Taylor expansion \eqref{Eq_Taylor_estimate_1} for every $x\in\cN_0$ into
\begin{align}
\Big|h_k(x)-\frac{1}{2}x^\top\cH\{h\}(0)x\Big|&=\frac{1}{2}\big|x^\top\big(\cH\{h_k\}(\xi_{k,x})-\cH\{h\}(0)\big)x\big| \notag \\
&\leq\frac{\big\Vert\cH\{h_k\}(\xi_{k,x})-\cH\{h\}(0)\big\Vert}{2h_0}\,x^\top\cH\{h\}(0)x. \label{Eq_Taylor_estimate_2}
\end{align}
Let us now choose $\varepsilon>0$ arbitrary. Since $h\in C^2(\cN_0)$, there exists an open neighborhood $\cN_\varepsilon$ of $x_0$, such that
\begin{equation*}
\Vert\cH\{h\}(\xi)-\cH\{h\}(0)\Vert<\varepsilon h_0,\qquad\xi\in\cN_\varepsilon.
\end{equation*}
Moreover, since there converges $\lim_{k\to\infty}\Vert h_k-h\Vert_{C^2(\cN_0)}$ by \Cref{ass_Laplace_assumption_discrete}~i), there exists some $K_\varepsilon\in\mathbb{N}$, such that
\begin{equation*}
\sup\limits_{\xi\in\cN_0}\Vert\cH\{h_k\}(\xi)-\cH\{h\}(\xi)\Vert<\varepsilon h_0,\qquad k\geq K_\varepsilon.
\end{equation*}
Combining these two estimates, we then get
\begin{align*}
\Vert\cH\{h_k\}(\xi_{k,x})-\cH\{h\}(0)\Vert&\leq\Vert\cH\{h_k\}(\xi_{k,x})-\cH\{h\}(\xi_{k,x})\Vert \\
&\quad+\Vert\cH\{h\}(\xi_{k,x})-\cH\{h\}(0)\Vert\leq 2\varepsilon h_0,\qquad x\in\cN_\varepsilon,\,k\geq K_\varepsilon.
\end{align*}
Using this in the estimate \eqref{Eq_Taylor_estimate_2}, then gives
\begin{equation*}
\Big|h_k(x)-\frac{1}{2}x^\top\cH\{h\}(0)x\Big|\leq\varepsilon\,x^\top\cH\{h\}(0)x,\qquad x\in\cN_\varepsilon,\,k\geq K_\varepsilon. \qedhere
\end{equation*}
\end{proof}

Now we are ready prove the convergence of the discrete Laplace method \Cref{thm_Discrete_Laplace_method}. It is a generalization of \cite[Lemma 3.6]{dBM-FOCS}, which they loosely justified as being provable by an
approach similar to that for
the continuous analog, e.g. \Cref{lemma:laplace-2}.

\begin{theorem}[Discrete Laplace method]\label{thm_Discrete_Laplace_method}
Let $\cW\subseteq\mathbb{R}^d$ be bounded, and for every $k\in\mathbb{N}$ consider a lattice $\Lambda_k$ as in Eq.\eqref{Eq_Lattice}, as well as functions $g_k:\cW\cap\Lambda_k\to\mathbb{C}$ and $h_k:\cW\to\mathbb{R}$, which satisfy Assumption~\ref{ass_Laplace_assumption_discrete}. Then there converges
\begin{equation}\label{Eq_Discrete_Laplace_method}
\lim\limits_{k\to\infty}n_k^{-\frac{d}{2}}e^{n_kh(x_0)}\sum_{x\in\cW\cap\Lambda_k}g_k(x)e^{-n_kh_k(x)}=\frac{(2\pi)^{\frac{d}{2}}g(x_0)}{|\det(A)|\sqrt{\det(\cH\{h\}(x_0))}}.
\end{equation}
\end{theorem}

\begin{proof}
First we note that for the (uniform) convergence of {a sequence of} complex functions, it is equivalent to ask for the (uniform) convergence of the respective positive and negative part of the real and complex part. Using this fact in the limit \eqref{Eq_Discrete_Laplace_method} but also in the Assumption~\ref{ass_Laplace_assumption_discrete}~i), we can without loss of generality assume that $g,g_k\geq 0$. For an easier notation let us also set $x_0=0$, $h(x_0)=0$ and use the abbreviation $H:=\cH\{h\}(x_0)$.

Let us fix $\varepsilon>0$ arbitrary. Then by Lemma \ref{lem_Taylor_estimate}, there exists some neighborhood $\cN_\varepsilon$ of $x_0$ and some $K_\varepsilon\in\mathbb{N}$, such that
\begin{equation}\label{Eq_Discrete_Laplace_method_1}
\frac{1-\varepsilon}{2}x^\top Hx\leq h_k(x)\leq\frac{1+\varepsilon}{2}x^\top Hx,\qquad x\in\cN_\varepsilon,\,k\geq K_\varepsilon.
\end{equation}
Since by \Cref{ass_Laplace_assumption_discrete}~i), there converges $\lim_{k\to\infty}\Vert g_k-g\Vert_{\infty,\,\cN_0\cap\Lambda_k}=0$ uniformly to a continuous function $g$, we can choose $\cN_\varepsilon\cN_0$ small enough and $K_\varepsilon$ large enough, such that also
\begin{equation*}
|g_k(x)-g(0)|\leq\Vert g_k-g\Vert_{\infty,\mathcal{N}_0\cap\Lambda_k}+|g(x)-g(0)|\leq\varepsilon,\qquad  x\in\cN_\varepsilon\cap\Lambda_k,\,k\geq K_\varepsilon.
\end{equation*}
This is equivalent to
\begin{equation}\label{Eq_Discrete_Laplace_method_2}
g(0)-\varepsilon\leq g_k(x)\leq g(0)+\varepsilon,\qquad  x\in\cN_\varepsilon\cap\Lambda_k,\,k\geq K_\varepsilon.
\end{equation}
With this preparation, we turn to estimating the sum. To do so, we split it up into the two terms
\begin{equation}\label{Eq_Discrete_Laplace_method_3}
\sum\limits_{x\in\cW\cap\Lambda_k}g_k(x)e^{-n_kh_k(x)}=\sum\limits_{x\in(\cW\setminus\cN_\varepsilon)\cap\Lambda_k}g_k(x)e^{-n_kh_k(x)}+\sum\limits_{x\in\cN_\varepsilon\cap\Lambda_k}g_k(x)e^{-n_kh_k(x)}.
\end{equation}
For the first sum, it is already shown in Lemma~\ref{lem_Sum_outside}, that there converges
\begin{equation}\label{Eq_Discrete_Laplace_method_4}
\lim\limits_{k\to\infty}n_k^{-\frac{d}{2}}\sum_{x\in(\cW\setminus\cN_\varepsilon)\cap\Lambda_k}g_k(x)e^{-n_kh_k(x)}=0.
\end{equation}
For the second sum in Eq.\eqref{Eq_Discrete_Laplace_method_3}, we use the first inequality in Eq.\eqref{Eq_Discrete_Laplace_method_1} and the second inequality in Eq.\eqref{Eq_Discrete_Laplace_method_2}, to estimate
\begin{equation*}
\sum_{x\in\cN_\varepsilon\cap\Lambda_k}g_k(x)e^{-n_kh_k(x)}\leq(g(0)+\varepsilon)\sum_{x\in\cN_\varepsilon\cap\Lambda_k}e^{-n_k\frac{(1-\varepsilon)}{2}x^\top Hx},\qquad k\geq K_\varepsilon.
\end{equation*}
If multiplying $n_k^{-\frac{d}{2}}$ and sending $k\to\infty$, \Cref{lem_Sum_Gaussian} tells us that the right hand side converges to
\begin{equation}\label{Eq_Discrete_Laplace_method_5}
\limsup\limits_{k\to\infty}n_k^{-\frac{d}{2}}\sum_{x\in\cN_\varepsilon\cap\Lambda_k}g_k(x)e^{-n_kh_k(x)}\leq\frac{(g(0)+\varepsilon)(2\pi)^{\frac{d}{2}}}{|\det(A)|\sqrt{(1-\varepsilon)^d\det(H)}}.
\end{equation}
Adding now the two limits \eqref{Eq_Discrete_Laplace_method_4} and \eqref{Eq_Discrete_Laplace_method_5}, gives the upper bound of the limit superior
\begin{equation}\label{Eq_Discrete_Laplace_method_6}
\limsup\limits_{k\to\infty}n_k^{-\frac{d}{2}}\sum_{x\in\cW\cap\Lambda_k}g_k(x)e^{-n_kh_k(x)}\leq\frac{(g(0)+\varepsilon)(2\pi)^{\frac{d}{2}}}{|\det(A)|\sqrt{(1-\varepsilon)^d\det(H)}}.
\end{equation}
Since this is true for every $\varepsilon>0$, which only appears on the right hand side, we can send $\varepsilon\to 0^+$, and conclude
\begin{equation}\label{Eq_Discrete_Laplace_method_7}
\limsup\limits_{k\to\infty}n_k^{-\frac{d}{2}}\sum_{x\in\cW\cap\Lambda_k}g_k(x)e^{-n_kh_k(x)}\leq\frac{g(0)(2\pi)^{\frac{d}{2}}}{|\det(A)|\sqrt{\det(H)}}.
\end{equation}
Analogously to \eqref{Eq_Discrete_Laplace_method_6}, we also derive the lower bound of the limit inferior
\begin{equation*}
\liminf\limits_{k\to\infty}n_k^{-\frac{d}{2}}\sum_{x\in\cW\cap\Lambda_k}g_k(x)e^{-n_kh_k(x)}\geq\frac{(g(0)-\varepsilon)(2\pi)^{\frac{d}{2}}}{|\det(A)|\sqrt{(1+\varepsilon)^d\det(H)}}.
\end{equation*}
Sending $\varepsilon\to 0^+$ also here, gives
\begin{equation}\label{Eq_Discrete_Laplace_method_8}
\liminf\limits_{k\to\infty}n_k^{-\frac{d}{2}}\sum_{x\in\cW\cap\Lambda_k}g_k(x)e^{-n_kh_k(x)}\geq\frac{g(0)(2\pi)^{\frac{d}{2}}}{|\det(A)|\sqrt{\det(H)}}.
\end{equation}
Since the upper bound of the limit superior \eqref{Eq_Discrete_Laplace_method_7} coincides with the lower bound of the limit inferior \eqref{Eq_Discrete_Laplace_method_8}, we conclude the limit \eqref{Eq_Discrete_Laplace_method}.
\end{proof}

\subsection{Asymptotics for discrete summands}\label{sec_Discrete_summands}

This section now generalizes the results from \Cref{thm_Discrete_Laplace_method} in the sense that the convergence of the sum \eqref{Eq_Discrete_Laplace_method} will be proven, not for the summands $g_k(x)e^{-n_kh_k(x)}$, but for discrete summands which in some sense only asymptotically behave as
\begin{equation*}
\cS_k(x)\sim g_k(x)e^{-n_kh_k(x)}.
\end{equation*}

\begin{theorem}[Discrete Laplace method on discrete summands]\label{thm_Discrete_Sk}
Let $\cW\subseteq\mathbb{R}^d$ be bounded, and $(\Lambda_k)_k$ a sequence of lattices as in Eq.\eqref{Eq_Lattice}. Let now $(\cS_k)_k:\cW\cap\Lambda_k\to\mathbb{C}$ be a sequence of functions, such that there exist for every $k\in\mathbb{N}$ functions $g_k:\cW\cap\Lambda_k\to\mathbb{C}$ and $h_k:\cW\to\mathbb{R}$, which satisfy the \Cref{ass_Laplace_assumption_discrete}, as well as the following two assumptions:
\index{$\cS_k(x)$} \medskip

\begin{enumerate}
\item[i)] There exists a constant $\mathdf{C_{\cS}}\geq 0$, such that
\begin{equation}\label{Eq_Sk_estimate1}
|\cS_k(x)|\leq C_{\cS}|g_k(x)|e^{-n_kh_k(x)},\qquad x\in\cW\cap\Lambda_k,\,k\in\mathbb{N}.
\end{equation}

\item[ii)] In the open neighborhood $\cN_0$ from \Cref{ass_Laplace_assumption_discrete}, there converges
\begin{equation}\label{Eq_Sk_estimate2}
\lim\limits_{k\to\infty}\Vert\cS_ke^{n_kh_k}-g_k\Vert_{\infty,\,\cN_0\cap\Lambda_k}=0.
\end{equation}
\end{enumerate}

Then there converges
\begin{equation}\label{Eq_Discrete_Sk}
\lim\limits_{k\to\infty}n_k^{-\frac{d}{2}}e^{n_kh(x_0)}\sum\limits_{x\in\cW\cap\Lambda_k}\cS_k(x)=\frac{(2\pi)^{\frac{d}{2}}g(x_0)}{\abs{\det(A)}\sqrt{\det(\mathcal{H}\{h\}(x_0))}}.
\end{equation}
\end{theorem}

\begin{proof}
Taking into account the convergence \eqref{Eq_Discrete_Laplace_method} of \Cref{thm_Discrete_Laplace_method}, it is enough to prove that
\begin{equation}\label{Eq_Discrete_Sk_1}
\lim\limits_{k\to\infty}n_k^{-\frac{d}{2}}e^{n_kh(x_0)}\sum\limits_{x\in\cW\cap\Lambda_k}|\cS_k(x)-g_k(x)e^{-n_kh_k(x)}|=0.
\end{equation}
Let us fix $\varepsilon>0$ arbitrary. Then by \Cref{lem_Taylor_estimate}, there exists an open neighborhood $\cN_\varepsilon$ of $x_0$ and some $K_\varepsilon\in\mathbb{N}$, such that the inequality \eqref{Eq_Taylor_estimate} holds true. Now, let us first consider the summation over $(\cW\setminus\cN_\varepsilon)\cap\Lambda_k$. With the estimate \eqref{Eq_Sk_estimate1}, we get
\begin{align*}
\sum\limits_{x\in(\cW\setminus\cN_\varepsilon)\cap\Lambda_k}|\cS_k(x)-g_k(x)e^{-n_kh_k(x)}|&\leq\sum\limits_{x\in(\cW\setminus\cN_\varepsilon)\cap\Lambda_k}\Big(|\cS_k(x)|+|g_k(x)|e^{-n_kh_k(x)}\Big) \\
&\leq(C_{\cS}+1)\sum\limits_{x\in(\cW\setminus\cN_\varepsilon)\cap\Lambda_k}|g_k(x)|e^{-n_kh_k(x)}.
\end{align*}
From Lemma \ref{lem_Sum_outside} we then conclude the convergence
\begin{align}
\lim\limits_{k\to\infty}n_k^{-\frac{d}{2}}&e^{n_kh(x_0)}\sum\limits_{x\in(\cW\setminus\cN_\varepsilon)\cap\Lambda_k}\big|\cS_k(x)-g_k(x)e^{-n_kh_k(x)}\big| \notag \\
&\leq(C_{\cS}+1)\lim\limits_{k\to\infty}n_k^{-\frac{d}{2}}e^{n_kh(x_0)}\sum\limits_{x\in(\cW\setminus\cN_\varepsilon)\cap\Lambda_k}|g_k(x)|e^{-n_kh_k(x)}=0. \label{Eq_Discrete_Sk_3}
\end{align}
For the sum over $\cN_\varepsilon\cap\Lambda_k$ in Eq.\eqref{Eq_Discrete_Sk_1}, we use the assumption \eqref{Eq_Sk_estimate2}, to possibly increase $K_\varepsilon$ and obtain
\begin{equation*}
\Vert\cS_ke^{n_kh_k}-g_k\Vert_{\infty,\,\cN_0\cap\Lambda_k}<\varepsilon,\qquad k\geq K_\varepsilon.
\end{equation*}
Use also the upper bound \eqref{Eq_Taylor_estimate} on the function $h_k$, and assume that $\cN_\varepsilon\subseteq\cN_0$, we obtain
\begin{align*}
\sum\limits_{x\in\cN_\varepsilon\cap\Lambda_k}\big|\cS_k(x)-g_k(x)&e^{-n_kh_k(x)}\big|=\sum\limits_{x\in\cN_\varepsilon\cap\Lambda_k}\big|\cS_k(x)e^{n_kh_k(x)}-g_k(x)\big|e^{-n_kh_k(x)} \\
&\leq\varepsilon\sum\limits_{x\in\cN_\varepsilon\cap\Lambda_k}e^{-n_k\big(h(x_0)+\frac{1-\varepsilon}{2}(x-x_0)^\top\cH\{h\}(x_0)(x-x_0)\big)},\qquad k\geq K_\varepsilon.
\end{align*}
If we now multiply $n_k^{-\frac{d}{2}}e^{n_kh(x_0)}$ and send $k\to\infty$, Lemma~\ref{lem_Sum_Gaussian} yields the convergence
\begin{align}
\limsup\limits_{k\to\infty}n_k^{-\frac{d}{2}}e^{n_kh(x_0)}\sum\limits_{x\in\cN_\varepsilon\cap\Lambda_k}&\big|\cS_k(x)-g_k(x)e^{-n_kh_k(x)}\big| \notag \\
&\leq\varepsilon\lim\limits_{k\to\infty}n_k^{-\frac{d}{2}}\sum\limits_{x\in\cN_\varepsilon\cap\Lambda_k}e^{-n_k\frac{1-\varepsilon}{2}(x-x_0)^\top\cH\{h\}(x_0)(x-x_0)} \notag \\
&=\frac{\varepsilon(2\pi)^{\frac{d}{2}}}{|\det(A)|\sqrt{(1-\varepsilon)^d\det(\cH\{h\}(x_0))}}. \label{Eq_Discrete_Sk_4}
\end{align}
Combining the limits \eqref{Eq_Discrete_Sk_3} and \eqref{Eq_Discrete_Sk_4} then gives the limit superior of the whole sum
\begin{equation*}
\limsup\limits_{k\to\infty}n_k^{-\frac{d}{2}}e^{n_kh(x_0)}\sum\limits_{x\in\cW\cap\Lambda_k}\big|\cS_k(x)-g_k(x)e^{-n_kh_k(x)}\big|\leq\frac{\varepsilon(2\pi)^{\frac{d}{2}}}{\abs{\det(A)}\sqrt{(1-\varepsilon)^d\det(\cH\{h\}(x_0))}}.
\end{equation*}
Since $\varepsilon>0$ was arbitrary, and the left hand side does not depend on $\varepsilon$, we can send $\varepsilon\to 0^+$ and obtain the desired vanishing limit \eqref{Eq_Discrete_Sk_1}.
\end{proof}

\section{Application: 3XOR-SAT}\label{sec_Application_3XORSAT_3XORGAME}

This section is devoted to giving an idea of how our theory, in particular \Cref{thm_Discrete_Sk}, can be applied. We describe two examples called 3XOR-SAT and 3XOR-GAME. Since the formulas of 3XOR-GAME are of similar structure but much more complicated, we will only consider the 3XOR-SAT problem in detail.

\subsection{3XOR solvability problems}

The 3XOR problems we shall discuss are equivalent to solving (binary) linear equations of the form
\begin{equation*}
\mathdf{\Gamma}z=b\pmod{2},
\end{equation*}
with $\Gamma$ a $m\times n$ matrix whose entries are all nonnegative integers, and $b$ is vector whose entries are all 0 or 1. In addition, the matrix $\Gamma$ has very special structure. \medskip

\textbf{3XOR-SAT equations} are those where $\Gamma$ is any XOR matrix, which means each row adds exactly to 3. Note that this is only possible for matrices $\Gamma$ with entries 0, 1, 2, or 3. \medskip

\textbf{3XOR-GAME equations} are those where the columns of $\Gamma$ partition into 3 parts 
\begin{equation*}
\Gamma=\begin{pmatrix} A & B & C \end{pmatrix},
\end{equation*}
with $A,B,C$ are $m\times\frac{n}{3}$ matrices whose rows each have exactly one 1, the rest are 0s. Clearly, every 3XOR-GAME equation is a 3XOR-SAT equation. \medskip

A matrix $\Gamma$ is called a \textit{2-core}, provided each column has a sum which is at least 2. If $\Gamma$ is not a 2-core, we can remove the columns with all 0s or only one 1 as well as remove appropriate rows, to obtain a linear system whose solvability is the same as the original system. With this motivation, we will focus efforts now on matrices $\Gamma$ which are 2-cores. \medskip

The interpretation of these systems as games will not be described in this article, since our only interest is in the solvability of the respective linear equations. However, without definitions we say that 3XOR games are 3 player (Alice, Bob, Charlie) cooperative games. Binary solutions to such a game's equation produces a perfect strategy for that game. Let us also mention that kXOR equations are defined similarly, but only $k=3$ is used here. \medskip

A major issue is the so called \textit{solvability threshold}, i.e. fix the size $m\times n$ of the system, and understand the solvability of randomly (uniformly distributed) generated 3XOR matrices $\Gamma$ and vectors $b$. In particular, we study the probability
\begin{equation}\label{Eq_Pmn}
\mathdf{P_{m,n}}:=P(\Gamma z=b\;\text{has a binary vector solution }z).
\end{equation}
A basic phenomena is the existence of a ``solvability threshold'' $\mathdf{c^*}$, which is the unique non negative number, such that as $m,n\to\infty$ and $\frac{m}{n}\to c$, the probability \eqref{Eq_Pmn} of randomly generated
$m\times n$ 3XOR problems converges to
\begin{equation}\label{Eq_cstar}
\lim\limits_{m/n\rightarrow c}P_{m,n}=\begin{cases} 1, & \text{if }c<c^*, \\ 0, & \text{if }c>c^*. \end{cases}
\end{equation}

Our interest in games came from a recent advance \cite{watts20203xor,WHZ2022satisfiability} on 3XOR-GAMEs showing real (as opposed to binary) solutions give a ``perfect quantum strategy'' for a 3XOR-GAME. Unfortunately, the classical game threshold had not been analyzed and doing so led to needing the discrete Laplace method.

\subsection{Determining the solvability threshold of 3XOR-SAT}

For the 3XOR-SAT problem, the existence and the value of the solvability threshold $\mathdf{c_{\text{SAT}}^*}$ (without the 2-core assumption on the matrix $\Gamma$) was shown in the classic paper \cite[Theorem 1.1]{dBM-FOCS}. This was a lovely tour de force with a many part argument, one of which was calculating the solvability threshold $\mathdf{c_{\text{SAT,2core}}^*}$ (with the 2-core assumption on the matrix $\Gamma$) in \cite[Theorem 3.1]{dBM-FOCS}, using the discrete Laplace method. The justification given was ``The proof is broadly similar to the ordinary Laplace approximation, (see, e.g., \cite{deBr81}), except that instead of working directly on the integral, one has to highlight a Riemann sum". See also Remark \ref{rem:csat} for more details. In deriving our critical threshold analog for 3XOR-GAMES we wanted a full proof; which is what motivated this paper, together with the hope that a general proved theorem might be more broadly useful.

The key to $c_{\text{SAT}}^*$ in \cite{dBM-FOCS} is: For fixed size $m \times n$ randomly generated 3XOR-SAT 2-core equations, it is shown in \cite[Lemma 3.2, Lemma 3.3]{dBM-FOCS} that the probability \eqref{Eq_Pmn} that a solution exists is bounded from below and from above by
\begin{equation}\label{Eq_P_bounds}
\bigg(\frac{1}{n}\sum_{(r,\alpha)\in\hW\cap\Lambda_{m,n}}\mathcal{S}_{m,n}(r,\alpha)\bigg)^{-1}\leq P_{m,n}\leq 2^{n-m},
\end{equation}
where the summand\footnote{For those looking at \cite{dBM-FOCS}, they omit this full summand, but give $\#\mathcal{F}\mathbb{E}[N^2]$ in pursuit of understanding the ratio $\frac{\mathbb{E}[N^2]}{\mathbb{E}[N]^2}$. Hence Eq.\eqref{Eq_Smn} arises by introducing $\#\mathcal{F}=2^mS_2(3m,n)n!$ and $\mathbb{E}[N]=2^{n-m}$. We also artificially introduce a factor $n$ in Eq.\eqref{Eq_Smn}, to cancel the pre-factor $\frac{1}{n}$ in Eq.\eqref{Eq_P_bounds}} is
\begin{equation}\label{Eq_Smn}
\mathdf{\mathcal{S}_{m,n}(r,\alpha)}:=n2^{m-n}3^{\frac{3}{2}(1-r)m}\binom{m}{\frac{3}{2}(1-r)m}\frac{S_2(3rm,\alpha n)S_2\big(3(1-r)m,(1-\alpha)n\big)}{S_2(3m,n)},
\end{equation}
where $S_2(p,q)$ denote the $2$-associated Stirling numbers of the second kind, and the summation region depends on the set $\hW$ as well as the lattice $\Lambda_{m,n}$, defined\footnote{Dubois-Mandler in \cite{dBM-FOCS} is slightly vague in their lattice definition and simply state $r\in\frac{1}{3m}\mathbb{Z}$. However, the detailed choice is $r\in\frac{2}{3m}\mathbb{Z}+1$, since $\frac{3}{2}(1-r)m$ in Eq.\eqref{Eq_Smn} must be an integer. This is because it arises as an integer counting the number of rows in a submatrix.} as
\begin{equation*}
\mathdf{\hW}=[\mathsmaller{\frac{1}{3}},1]\times[0,1]\qquad\text{and}\qquad\mathdf{\Lambda_{m,n}}:=\begin{pmatrix} \frac{2}{3m} & 0 \\ 0 & \frac{1}{n} \end{pmatrix}\mathbb{Z}^2+\begin{pmatrix} 1 \\ 0 \end{pmatrix}.
\end{equation*}
We now state the main theorem of this section.

\begin{theorem}[Application to $2$-core 3XOR-SAT]
\label{thm_Application_3XORSAT}
Let $(m_k)_k$ and $(n_k)_k$ be arbitrary sequences, with $m_k,n_k\geq k$ and $\lim_{k\to\infty}\frac{m_k}{n_k}\in(\frac{2}{3},1)$. Then
\begin{equation}\label{Eq_Application_3XORSAT}
\lim\limits_{k\to\infty}\frac{1}{n_k}\sum\limits_{(r,\alpha)\in\hW\cap\Lambda_{m_k,n_k}}\mathcal{S}_{m_k,n_k}(r, \alpha)=1.
\end{equation}
\end{theorem}

\begin{proof}
A sketch of the proof is given in Section~\ref{sec_Sketch_of_the_proof}. A more complete proof is given in the \online \Cref{sec_Appendix}. 
\end{proof}

Let us now return to the solvability threshold of the 3XOR-SAT problem. More precisely, for the 3XOR-SAT problem there exist two separate thresholds, $c_{\text{SAT}}^*$ corresponds the all possible $\Gamma$-matrices, and $c_{\text{SAT,2core}}^*$ is restricted to 2-core matrices $\Gamma$. In the upcoming \Cref{cor:Threshold}, and in the rest of this section, we shall focus on the 2-core threshold only, while the general threshold is mentioned in \Cref{rem:csat}. In accordance with the general definition of the critical threshold in Eq.\eqref{Eq_cstar}, the 2-core threshold is the unique number $c_{\text{SAT,2core}}^*\in(\frac{2}{3},\infty)$, such that as $m,n\to\infty$ and $\frac{m}{n}\to c$, the probability \eqref{Eq_Pmn} with randomly generated 2-core matrices $\Gamma$ behaves as
\begin{equation*}
\lim\limits_{m/n\rightarrow c}P_{m,n}=\begin{cases} 1, & \text{if }\frac{2}{3}<c<c_{\text{SAT,2core}}^*, \\ 0, & \text{if }c>c_{\text{SAT,2core}}^*. \end{cases}
\end{equation*}
Note the matrix $\Gamma$ being a 2-core forces $\frac{m}{n}\geq\frac{2}{3}$.

\begin{corollary}\label{cor:Threshold}
The solvability threshold for the $2$-core 3XOR-SAT problem is given by
\begin{equation*}
c_{\textnormal{SAT,2core}}^*=1.
\end{equation*}
\end{corollary}

\begin{proof}
Consider two sequences $(m_k)_k$ and $(n_k)_k$ of number of equations and number of variables, with $m_k,n_k\geq k$ and $\lim_{k\to\infty}\frac{m_k}{n_k}=:c$. \medskip

If $c\in(\frac{2}{3},1)$, it follows from the first inequality in \eqref{Eq_P_bounds} and from Theorem~\ref{thm_Application_3XORSAT}, that
\begin{equation*}
\lim\limits_{k\to\infty}P_{m_k,n_k}\geq\lim\limits_{k\to\infty}\bigg(\frac{1}{n_k}\sum_{(r,\alpha)\in\hW\cap\Lambda_{m_k,n_k}}\mathcal{S}_{m_k,n_k}(r,\alpha)\bigg)^{-1}=1.
\end{equation*}
If $c>1$ on the other hand, it follows from the second inequality in \eqref{Eq_P_bounds}, that
\begin{equation*}
\lim\limits_{k\to\infty}P_{m_k,n_k}\leq\lim\limits_{k\to\infty}2^{n_k-m_k}=\lim\limits_{k\to\infty}2^{-n_k(\frac{m_k}{n_k}-1)}=\lim\limits_{k\to\infty}2^{-n_k(c-1)}=0.
\end{equation*}
These two limits now show that $\lim_{k\to\infty}P_{m_k,n_k}=1$ if $c\in(\frac{2}{3},1)$ and $\lim_{k\to\infty}P_{m_k,n_k}=0$ if $c>1$, which proves $c_{\text{SAT,2core}}^*=1$ by the definition of the critical threshold in \eqref{Eq_cstar}.
\end{proof}

\begin{remark}\label{rem:csat}
Although the threshold $c_{\text{SAT}}^*$ will not be computed in this article, let us shortly explain its two key parts. One is \Cref{cor:Threshold}, which handles 2-cores, while the second uses a long separate argument, to analyze an algorithm which basically deletes certain rows and columns from a given matrix $\Gamma$, to end up with the (unique)
maximal 2-core submatrix $\widetilde \Gamma$.
If $\Gamma$ is uniformly distributed,
then so is $\widetilde{\Gamma}$; this principle is called ``maintenance of uniformity". The idea in applying these, c.f. \cite{dBM-FOCS} or 
\cite{HHgameThreshold}, is to select 
$0 < c < c^*_{\text{SAT}}$, and pick very large values $m,n$ with $\frac{m}{n}<c$. Corresponding to $c$, there then exists a value $\widetilde{c}<c^*_{\text{SAT,2core}}=1$, such that for any randomly generated $m \times n$ matrix $\Gamma$ the probability that its maximal 2-core
$\widetilde\Gamma$, has size $\widetilde m \times \widetilde n$ with $\frac{2}{3}\leq\frac{\widetilde m}{\widetilde n}<\widetilde c$, is very high. Hence \Cref{cor:Threshold} implies that there is a very high probability that a randomly generated $m \times n$ 3XOR-SAT problem has a solution. A similar argument settles the 
$c>c^*_{\text{SAT}}$ case.
\end{remark}

\subsection{Sketch of the proof of \Cref{thm_Application_3XORSAT}}\label{sec_Sketch_of_the_proof}

As we soon see, our
\Cref{thm_Discrete_Sk} on such sums 
concludes with the desired answer.
We describe this outcome in
\Cref{thm_Application_3XORSAT},
and to do this, we introduce some definitions.
First we discuss the asymptotic expansion
\begin{equation}\label{Eq_Stirling_type_approximation}
\cS_{m,n}(r,\alpha)\sim g_{m,n}(r,\alpha)e^{-nh_{m/n}(r,\alpha)},
\end{equation}
as $m,n\to\infty$ with $\frac{m}{n}\rightarrow c\in(\frac{2}{3},1)$. For the most part, the definitions of $g_{m,n}$ and $h_{m/n}$ are based directly on functions defined in \cite{dBM-FOCS}. The formulas from \cite{dBM-FOCS} are properly defined in the interior of the trapezoid
\begin{align*}
\mathdf{\hT_{3c}}\coloneqq \Big\{(r,\alpha) \in\widehat{ \mathcal{W}}\;\Big|\; 1 - \frac{3(1-r)c}{2} \leq \alpha\leq\frac{3rc}{2} \Big\}.
\end{align*}
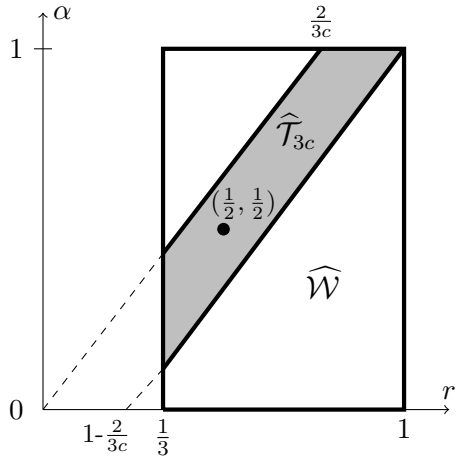
\begin{figure}[H]
\centering
\begin{tikzpicture}[scale=1.2]
\fill[black!25] (4,4)--(3.08,4)--(1.33,1.73)--(1.33,0.45);
\draw[->] (0,0)--(4.5,0) node[anchor=south] {$r$};
\draw[->] (0,0)--(0,4.4) node[anchor=west] {$\alpha$};
\draw[ultra thick] (1.33,0)--(4,0)--(4,4)--(1.33,4)--(1.33,0);
\draw (1.33,0) node[anchor=north] {$\frac{1}{3}$};
\draw (4,0) node[anchor=north] {$1$};
\draw (-0.1,0) node[anchor=east] {$0$};
\draw (0.1,4)--(-0.1,4) node[anchor=east] {$1$};
\draw[dashed] (0,0)--(1.33,1.73);
\draw[ultra thick] (1.33,1.73)--(3.08,4) node[anchor=south] {\small{$\frac{2}{3c}$}};
\draw[dashed] (0.92,0)--(1.33,0.45);
\draw[ultra thick] (1.33,0.45)--(4,4);
\draw (0.7,0) node[anchor=north] {\small{$1$-$\frac{2}{3c}$}};
\fill[black] (2,2) circle (0.07cm);
\draw (2.22,2) node[anchor=south] {\small{$(\frac{1}{2},\frac{1}{2})$}};
\draw (2.8,3.4) node[anchor=north] {\large{$\hT_{3c}$}};
\draw (3.1,1.7) node[anchor=north] {\large{$\widehat{\mathcal{W}}$}};
\end{tikzpicture}
\caption{Illustration of trapezoidal region $\hT_{3c}$, where the function $g_{m,n}(r,\alpha)$ is nonzero. Outside of the trapezoid $\hT_{3c}$, the function vanishes, $g_{m,n}(r,\alpha)=0$. Moreover, the interior of $\hT_{3c}$ contains the interesting part of $g_{m,n}(r,\alpha)$, where it also coincides with the definition in \cite{dBM-FOCS}. The definition of this function can also be found in Eq.\eqref{Eq_gmn} in the \online appendix of this article.}
\end{figure}

On the interior of the trapezoid $\hT_{3c}$, Dubois-Mandler represents $g_{m,n}$ and $h_{m/n}$ as
\begin{equation}\label{Eq:_gmn_hmn_Dubois_Mandler}
g_{m,n}=g_{m/n} \qquad \text{and} 
\qquad h_{m/n} = 2(1-\mathsmaller{\frac{m}{n}})\ln(2)-f_{m/n}.
\end{equation}
for some functions $g_c,f_c$ depending on only the one index $c = m/n$.

The confirmation of the hypotheses of 
\Cref{thm_Discrete_Sk}, a result which is needed to prove 
\Cref{thm_Application_3XORSAT}, is outlined
in two lemmas. A serious understanding of both lemmas 
would require the more detailed definitions
from the \online \Cref{sec_Appendix}.

\begin{lemma}[\cite{dBM-FOCS}]
\label{lem:duBMassumpIsTrue}
For arbitrary two sequences $(n_k)_k$, $(m_k)_k$, with $m_k,n_k\geq k$ and $\lim_{k\to\infty}\frac{m_k}{n_k}=c\in(\frac{2}{3},1)$, the functions $(g_{m_k,n_k})_k$ and $(h_{m_k/n_k})_k$ from Eq.\eqref{Eq:_gmn_hmn_Dubois_Mandler} satisfy the Assumption~\ref{ass_Laplace_assumption_discrete} at the point $(r_0,\alpha_0)=(\frac{1}{2},\frac{1}{2})$.
\end{lemma}

\begin{proof}
We list each item in 
\Cref{ass_Laplace_assumption_discrete}
and indicate the nature of the respective proof.
Proofs were given in \cite{dBM-FOCS} with
varying degrees of detail, though item iv) (because it involves the boundary of $\mathcal{T}_{3c}$)
requires different techniques; also the proof of our item iii) has some differences. This, and full definitions 
of functions $g_{m,n}, h_{m/n}$
are in \online \Cref{sec_Appendix}. \medskip

i)\;\;The convergence 
\begin{equation*}
\lim_{k\to\infty}\norm{g_{m_k,n_k}-g_c}_{\infty,\cN_0\cap\Lambda_{m_k,n_k}}=0\qquad\text{and}\qquad\lim_{k\to\infty}\norm{h_{m_k/n_k}-h_c}_{C^2(\cN_0)} = 0,
\end{equation*}
in a neighborhood $\mathcal{N}_0$ of $(r_0,\alpha_0)=(\frac{1}{2},\frac{1}{2})$, follows directly from the continuous $c$-dependence of $g_c$ and $h_c$. This can be seen directly from the explicit forms of the respective functions. \medskip

ii)\;\;Also the specific values $h_{m_k/n_k}(\frac{1}{2},\frac{1}{2})=h_c(\frac{1}{2},\frac{1}{2})=0$, and $\nabla h_{m_k/n_k}(\frac{1}{2},\frac{1}{2})=0$, as well as the positive definteness of the Hessian $\cH\{h_c\}(\frac{1}{2},\frac{1}{2})$ can be done by straight forward computation. \medskip

iii)\;\;For every open neighborhood $\mathcal{N}$ of $(r_0,\alpha_0)=(\frac{1}{2},\frac{1}{2})$, the lower bound
\begin{equation}\label{Eq_hc_minimum}
\inf_{k\in\mathbb{N}}\inf_{x\in\widehat{\cW}\setminus\cN}h_{m_k/n_k}(r,\alpha)>h_c(\mathsmaller{\frac{1}{2},\frac{1}{2}})=0,
\end{equation}
is ``obvious" from looking at the graph of $h_c$. Actually proving it is challenging, see \online \Cref{sec_Positivity_hy}. Note, that Eq.\eqref{Eq_hc_minimum} only holds for $c<1$, because $h_c(1,1)=(1-c)\ln(2)$ becomes negative otherwise, and $(\frac{1}{2},\frac{1}{2})$ is not the global minimum anymore. \medskip

iv)\;\;The fact that there exists a constant $C\geq 0$, such that
\begin{equation*}
\Vert g_{m_k,n_k}\Vert_{\infty,\widehat{\mathcal{W}}\cap\Lambda_{m_k,n_k}}\leq C n_k^3,\qquad k\in\mathbb{N},
\end{equation*}
follows laboriously from \Cref{lem_gS_estimate} and the calculations done in \online \Cref{sec_Subsection1}~iv).
\end{proof}

Now that we have shown that $g_{m_k,n_k},h_{m_k/n_k}$ satisfy \Cref{ass_Laplace_assumption_discrete}, the following lemma shows that the Stirling-type approximation \eqref{Eq_Stirling_type_approximation} is a good approximation to the summand $\cS_{m,n}$ in Eq.\eqref{Eq_Smn}. To be precise, we next show that the functions $\cS_{m_k,n_k}$ satisfy the hypotheses i) and ii) of \Cref{thm_Discrete_Sk}.

\begin{lemma}\label{lem:summandStirling}
Let $(n_k)_k$, $(m_k)_k$ be sequences with $m_k,n_k\geq k$ and $\lim_{k\to\infty}\frac{m_k}{n_k}=:c\in(\frac{2}{3},1)$.

\begin{enumerate} 
\item[i)] There exists a constant $C_{\cS}\geq 0$, such that for every $k\in\mathbb{N}$ there holds
\begin{equation*}
|\cS_{m_k,n_k}(r,\alpha)|\leq C_{\cS}|g_{m_k,n_k}(r,\alpha)|e^{-n_kh_{m_k/n_k}(r,\alpha)},\qquad (r,\alpha)\in\widehat{\cW}\cap\Lambda_{m_k,n_k}.
\end{equation*}

\item[ii)] There there exists an open neighborhood $\cN_0$ of $(r_0,\alpha_0)=(\frac{1}{2},\frac{1}{2})$, on which there converges
\begin{equation}\label{Eq_Sk_property2}
\lim\limits_{k\to\infty}\big\Vert\cS_{m_k,n_k}e^{n_kh_{m_k/n_k}}-g_{m_k,n_k}\big\Vert_{\infty,\,\cN_0\cap\Lambda_{m_k,n_k}}=0.
\end{equation}
\end{enumerate}
\end{lemma}

\begin{proof}
The approximation \eqref{Eq_Sk_property2} of $\cS_{m,n}(x)$ by $g_{m,n}(x) e^{-n h_{m/n}(x)}$ is derived by applying two asymptotics. One is the Stirling approximation $n! \sim \sqrt{2\pi n} (n/e)^n$, which converges as $n\to\infty$ and is correct within a constant factor for $n\geq 1$. The other is the Hennecart approximation for the Stirling numbers $S_2(p,q)$ put forth in \cite{Hennecart}, with uniform convergence except along the boundary of $\mathcal{T}_{3c}$, where $g_{m,n}$ blows up. The uniform convergence was conjectured in \cite{Hennecart} and assumed in \cite{dBM-FOCS}, but it was not proven until \cite{CHHarxiv}, a paper which relies heavily on \cite{stirlingConnDobro}. However, all approximations converge uniformly as needed in a neighborhood $\cN_0$ of $(r_0,\alpha_0)=(\half,\half)$, as needed to prove \Cref{Eq_Sk_property2}. 
\end{proof}

Now we are ready to prove the initial theorem of this section.

\begin{proof}[Proof of \Cref{thm_Application_3XORSAT}]
The idea is to apply \Cref{thm_Discrete_Sk}, to evaluate explicitly the value of the limit \eqref{Eq_Application_3XORSAT}. By \Cref{lem:duBMassumpIsTrue}, the \Cref{ass_Laplace_assumption_discrete} of the functions $(g_{m_k,n_k})_k$ and $(h_{m_k/n_k})_k$ is established with respect to the point $(r_0,\alpha_0)=(\half,\half)$. Also the assumptions i) and ii) of \Cref{thm_Discrete_Sk} of the summands $\mathcal{S}_{m_k,n_k}$ are satisfied due to \Cref{lem:summandStirling}.

The formula for the limit of the sum in Eq.\eqref{Eq_Discrete_Sk} then is
\begin{equation}\label{Eq_Final_limit}
\lim\limits_{k\to\infty}\frac{1}{n_k}\sum\limits_{ x\in\widehat{\cW}\cap\Lambda_{m_k,n_k}}\cS_{m_k,n_k}(x)=\frac{2\pi \, g_c(\half, \half )}{\abs{\det(A)}\sqrt{\det(\cH\{h_c\}(\half, \half))}}.
\end{equation}
The explicit values of the terms on the right hand side are known explicitly, see the algebra after Eq.(10) of \cite{dBM-FOCS}, or \online \Cref{sec_Subsection3} of our article. With these values, the right-hand side of Eq.\eqref{Eq_Final_limit} turns out to be equal to $1$, as claimed in the statement \eqref{Eq_Application_3XORSAT} of the theorem.
\end{proof}

\subsection{3XOR-GAME}\label{sec:3XOR-GAME}

To analyze solvability thresholds for 
3XOR-GAME equations, the first two authors of this article
tried the original approach, \cite{dBM-FOCS}.
Counting the number of relevant events plus
combinatorial arguments was used to compute key probabilities
and they yielded a sum
$\rho^{\text{GAME}}_k$ which is much more complicated than for SAT;
for example, it has 6 variables.
Very long arguments validate the
hypothesis of 
\Cref{thm_Discrete_Sk} when
$c<1$.
 As with 3XOR-SAT after much manipulation
 they got $\rho^{\text{GAME}}_k\sim 1 $ for 2-core games 
 and conclude that 3XOR-SAT and 3XOR-GAME
 have the same solvability threshold.

 This is a subject we do not belabor here,
 since the solution to the 3XOR-GAME solvability threshold, they plan to publish, is 
based on a different line of proof, see the arXiv-version \cite{HHgameThreshold}.

\newpage

\appendix

\section{Appendix}\label{sec_Appendix}

In the breakthrough paper \cite{dBM-FOCS}, some analytical details were described very briefly. Using \Cref{thm_Discrete_Sk} and recent theorems in \cite{stirlingConnDobro} and \cite{CHHarxiv}, these details 
can be thoroughly filled in. The goal here is to accomplish this,
and give full details of the proof of Theorem \ref{thm_Application_3XORSAT}. 
Since this appendix to the journal article is online only, our aim is to make its presentation fairly self contained, 
despite  entailing duplication with the main paper. For this reason we introduce the main objects this appendix will deal with. \medskip

Let us start with the \textbf{$2$-associated Stirling numbers of second kind}. For every $p,q\in\mathbb{N}_0$ we define the \textit{$2$-associated Stirling number $S_2(p,q)$} as the number of ways to partition a set of $p$ elements into $q$ subsets with each containing at least 2 elements. These numbers can inductively be written as
\begin{equation}\label{Eq_Stirling_number}
\begin{split}
{S_2(0,0)}&:=1;\qquad S_2(p,0):=0,\quad p\geq 1;\qquad S_2(p,q):=0,\quad p\geq 0,\,q>\frac{p}{2}; \\
\mathdf{S_2(p,q)}&:=qS_2(p-1,q)+(p-1)S_2(p-2,q-1),\quad 
{p\geq 2},\,1\leq q\leq\frac{p}{2}.
\end{split}
\end{equation}
If we now also consider the domain and the lattice
\begin{equation}\label{Eq_Lambdamn}
\hW:=[1/3,1]\times[0,1],\qquad\text{and}\qquad\Lambda_{m,n}:=\begin{pmatrix} \frac{2}{3m} & 0 \\ 0 & \frac{1}{n} \end{pmatrix}\mathbb{Z}^2+\begin{pmatrix} 1 \\ 0 \end{pmatrix},\qquad m,n\in\mathbb{N},
\end{equation}
we are able to define the main objectives of this appendix, namely for every $m,n\in\mathbb{N}$ the function $\mathcal{S}_{m,n}:\widehat{W}\cap\Lambda_{m,n}\to\mathbb{R}$, defined as
\begin{equation}\label{Eq_Smn_appendix}
\mathdf{\mathcal{S}_{m,n}(r,\alpha)}:=n2^{m-n}3^{\frac{3}{2}(1-r)m}\binom{m}{\frac{3}{2}(1-r)m}\frac{S_2(3rm,\alpha n)S_2\big(3(1-r)m,(1-\alpha)n\big)}{S_2(3m,n)}.
\end{equation}
With these objects, we are now in the situation to recall the statement of \Cref{thm_Application_3XORSAT}.

\begin{theorem}\label{thm_Application_3XORSAT_appendix}
Let $(m_k)_k$, $(n_k)_k$ be sequences, with $m_k,n_k\geq k$ and $\lim_{k\to\infty}\frac{m_k}{n_k}\in(\frac{2}{3},1)$. Then
\begin{equation}\label{Eq_Application_3XORSAT_appendix}
\lim\limits_{k\to\infty}\frac{1}{n_k}\sum\limits_{(r,\alpha)\in\hW\cap\Lambda_{m_k,n_k}}\mathcal{S}_{m_k,n_k}(r, \alpha)=1.
\end{equation}
\end{theorem}

This appendix will start with some standard notion of asymptotics in Section~\ref{sec_Asymptotics}, followed by the analysis of the individual terms of the summand $\mathcal{S}_{m,n}$ in Section~\ref{sec_Ingredients_of_Smn}, combining all this terms to get estimates and an asymptotic behaviour of the full summand $\mathcal{S}_{m,n}$ in Section~\ref{sec_Asymptotics_of_Smn}, proves positivity of the function $h_y$ in Section~\ref{sec_Positivity_hy}, and ends in Section~\ref{sec_Proof} with the actual proof of Theorem~\ref{thm_Application_3XORSAT_appendix}.

\subsection{Standard asymptotic notation}\label{sec_Asymptotics}

Here we will give a short review of standard notation for asymptotics.

\begin{definition}\label{defi_Asymptotics}
i)\;\;Let $(a_k)_k,(b_k)_k\in\mathbb{C}$ be sequences of complex numbers. We say that
\begin{equation*}
a_k\sim b_k,\qquad\text{as }k\to\infty,
\end{equation*}
if and only if
\begin{equation*}
\lim\limits_{k\to\infty}\Big|\frac{b_k}{a_k}-1\Big|=0.
\end{equation*}
ii)\;\;Let $a_k,b_k:\cW_k\rightarrow\mathbb{C}$ be two sequences of functions, defined on some set $\mathcal{W}_k$. We say that
\begin{equation}\label{Eq_Asymptotics_uniform}
a_k(x)\sim b_k(x),\qquad\text{as }k\to\infty,\text{ uniformly for }x\in\cW_k,
\end{equation}
if and only if
\begin{equation*}
\lim\limits_{k\to\infty}\sup\limits_{x\in\cW_k}\Big|\frac{b_k(x)}{a_k(x)}-1\Big|=0.
\end{equation*}
\end{definition}

The notion of asymptotics in Definition \ref{defi_Asymptotics} now has the following useful properties.

\begin{lemma}
Let $(a_k)_k,(b_k)_k,(c_k)_k,(d_k)_k\in\mathbb{C}$ are sequences. Then the notion $\sim$ of asymptotics in Definition \ref{defi_Asymptotics} has the following properties: \medskip

\begin{enumerate}
\item[i)] If $a_k\sim b_k$ and $b_k\sim c_k$, as $k\to\infty$, then also $a_k\sim c_k$, as $k\to\infty$. \medskip

\item[ii)] If $a_k\sim b_k$ and $c_k\sim d_k$, as $k\to\infty$, then also $a_kc_k\sim b_kd_k$, as $k\to\infty$. \medskip

\item[iii)] If $a_k\sim b_k$, as $k\to\infty$, then also $\frac{1}{a_k}\sim\frac{1}{b_k}$, as $k\to\infty$. \medskip
\end{enumerate}

If $a_k,b_k,c_k,d_k:\cW_k\rightarrow\mathbb{C}$ are sequences of functions, the above results i) and ii) are still valid if one considers $\sim$ as in Eq.\eqref{Eq_Asymptotics_uniform} uniformly for $x\in\cW_k$. The property iii) needs the additional assumption $\sup_{k\in\mathbb{N}}\sup_{x\in\cW_k}\frac{|a_k(x)|}{|b_k(x)|}<\infty$.
\end{lemma}

\subsection{Ingredients of the summands $\mathcal{S}_{m,n}(r,\alpha)$ and estimates on them}\label{sec_Ingredients_of_Smn}

A central role in this subsection will play the function
\begin{equation}\label{Eq_Q}
\mathdf{Q(x)}:=\frac{x(e^x-1)}{e^x-1-x},\qquad x>0.
\end{equation}
Since $Q$ is strictly increasing, with limits $\lim_{x\to 0^+}Q(x)=2$ and $\lim_{x\to\infty}Q(x)=\infty$, there exists the (again strictly increasing) inverse function
\index{$Q^{-1}(\xi)$}
\begin{equation*}
Q^{-1}\colon(2,\infty)\rightarrow(0,\infty).
\end{equation*}
With these functions $Q$ and $Q^{-1}$, and also their combinations
\begin{align}
\mathdf{P(\xi)}:=&\sqrt{Q^{-1}(\xi)Q'(Q^{-1}(\xi))},\qquad\xi\in(2,\infty), \label{Eq_P} \\
\mathdf{A(\xi)}:=&e^{Q^{-1}(\xi)}-1-Q^{-1}(\xi),\qquad\xi\in(2,\infty), \label{Eq_A}
\end{align}
let us define the functions
\begin{align}
\mathdf{g^S(p,q)}:=&\begin{cases} \frac{1}{\sqrt{2\pi q}\,P(\frac{p}{q})}, & 1\leq q<\frac{p}{2}, \\ 1, & q=\frac{p}{2}, \\ 0, & q=0\text{ or }q>\frac{p}{2}, \end{cases}\qquad 
p\in\mathbb{N},\,q\in\mathbb{N}_0,
\label{Eq_gS} \\
\mathdf{h^S(\xi)}:=&\begin{cases} \xi\ln(Q^{-1}(\xi))-\ln(A(\xi)), & \xi\in(2,\infty), \\ \ln(2), & \xi\in(0,2]\cup\{\infty\}. \end{cases} \label{Eq_hS}
\end{align}
Note that one can simply show that the function $h^S$ is continuous in the point $\xi=2$, but not at $\infty$, since there is $\lim_{\xi\rightarrow\infty}h^S(\xi)=\infty\neq\ln(2)=h^S(\infty)$. 

\begin{proposition}
i)\;\;There exist constants $0<c_S\leq C_S$, such that the Stirling numbers satisfy the estimate
\index{$c_s$, $C_S$}
\begin{equation}\label{Eq_S2_estimate}
c_S\frac{p!}{q!}g^S(p,q)e^{-qh^S(\frac{p}{q})}\leq S_2(p,q)\leq C_S\frac{p!}{q!}g^S(p,q)e^{-qh^S(\frac{p}{q})},\qquad p\in\mathbb{N},\,q\in\mathbb{N}_0.
\end{equation}
ii)\;\;Let $(m_k)_k,(n_k)_k\subseteq\mathbb{N}$, with $\lim_{k\to\infty}m_k=\lim_{k\to\infty}n_k=\infty$ and $\lim_{k\to\infty}\frac{m_k}{n_k}>\frac{2}{3}$. Then there exists a neighborhood $\mathcal{N}_0$ of $(\frac{1}{2},\frac{1}{2})$, such that
\begin{equation}\label{Eq_S2_approximation}
S_2(3rm_k,\alpha n_k)\sim\frac{(3rm_k)!}{(\alpha n_k)!}g^S(3rm_k,\alpha n_k)e^{-\alpha n_kh^S(\frac{3rm_k}{\alpha n_k})},\quad\text{as }k\to\infty,
\end{equation}
uniformly for $(r,\alpha)\in\mathcal{N}_0\cap\Lambda_{m_k,n_k}$. Here the asymptotic $\sim$ is understood in the sense of Eq.\eqref{Eq_Asymptotics_uniform}.
\end{proposition}

\begin{proof}
i)\;\;For $q=\frac{p}{2}$, one can easily verify inductively from Eq.\eqref{Eq_Stirling_number} the specific Stirling numbers $S_2(2q,q)=\frac{(2q)!}{q!2^q}$. With the explicit values $g^S(2q,q)=1$ and $h^S(\frac{2q}{q})=\ln(2)$, we then conclude the equality
\begin{equation*}
S_2(p,q)=\frac{(2q)!}{q!2^q}=\frac{p!}{q!}g^S(p,q)e^{-qh^S(\frac{p}{q})}.
\end{equation*}
For $q>\frac{p}{2}$, the Stirling numbers $S_2(p,q)=0$ vanish, and likewise for $q=0$ (since the case $p=0$ is not covered by the statement \eqref{Eq_S2_estimate}). In both cases also $g^S(p,q)=0$, and hence \eqref{Eq_S2_estimate} is satisfied with equality. For the crucial region $1\leq q<\frac{p}{2}$, we know from \cite[Theorem 1.2.1]{CHHarxiv} that there exists a function $E\colon\mathbb{N}\to\mathbb{R}$, with $\lim_{p\to\infty}E(p)=0$, which is an upper bound of
\begin{equation}\label{Eq_S2_asymptotics_3}
\Big|\frac{\mathdf{F(p,q)}}{S_2(p,q)}-1\Big|\leq\mathdf{E(p)},\qquad 1\leq q<\frac{p}{2},
\end{equation}
where the function $F$ is given by
\begin{equation*}
F(p,q):=\frac{p!\sqrt{2\pi(p-2q)}}{q!(p-2q)!}\Big(\frac{p-2q}{e}\Big)^{p-2q}g^S(p,q)e^{-qh^S(\frac{p}{q})}.
\end{equation*}
This estimate on the one hand shows that
\begin{equation}\label{Eq_S2_asymptotics_5}
\frac{F(p,q)}{S_2(p,q)}\leq\bigg|\frac{F(p,q)}{S_2(p,q)}-1\bigg|+1\leq E(p)+1\leq\sup\limits_{p\in\mathbb{N}}E(p)+1=:\frac{1}{c_F}.
\end{equation}
On the other hand, we can choose $N\in\mathbb{N}$ large enough, such that $E(p)\leq\frac{1}{2}$ for every $p\geq N$. For those values of $p$ we then get
\begin{equation*}
\frac{F(p,q)}{S_2(p,q)}\geq\frac{1}{2},\qquad p\geq N,\, 1\leq q<\frac{p}{2}.
\end{equation*}
If we now choose $\frac{1}{C_F}:=\min\Big\{\frac{1}{2}\,,\,\min_{1\leq q<\frac{p}{2}<\frac{N}{2}}\frac{F(p,q)}{S_2(p,q)}\Big\}>0$, we get the lower bound
\begin{equation}\label{Eq_S2_asymptotics_6}
\frac{F(p,q)}{S_2(p,q)}\geq\frac{1}{C_F},\qquad 1\leq q<\frac{p}{2}.
\end{equation}
Combining now Eq.\eqref{Eq_S2_asymptotics_5} and \eqref{Eq_S2_asymptotics_6} gives the upper and lower bound of the Stirling number
\index{$c_F$, $C_F$}
\begin{equation}\label{Eq_S2_asymptotics_7}
c_FF(p,q)\leq S_2(p,q)\leq C_FF(p,q).
\end{equation}
Next, from the Stirling approximation
\begin{equation}\label{Eq_S2_asymptotics_4}
\Big|\frac{n!}{\sqrt{2\pi n}}\Big(\frac{e}{n}\Big)^n-1\Big|\leq e^{\frac{1}{12n}}-1,\qquad n\in\mathbb{N},
\end{equation}
there follows the estimate
\begin{equation}\label{Eq_Factorial_estimate}
2-e^{\frac{1}{12}}\leq\frac{n!}{\sqrt{2\pi n}}\Big(\frac{e}{n}\Big)^n\leq e^{\frac{1}{12}},\qquad n\in\mathbb{N}.
\end{equation}
Using this inequality with $n=p-2q$, then gives
\begin{equation}\label{Eq_S2_asymptotics_8}
2-e^{\frac{1}{12}}\leq\frac{\frac{p!}{q!}g^S(p,q)e^{-qh^S(\frac{p}{q})}}{F(p,q)}=\frac{(p-2q)!}{\sqrt{2\pi(p-2q)}}\Big(\frac{e}{p-2q}\Big)^{p-2q}\leq e^{\frac{1}{12}}.
\end{equation}
Combining now Eq.\eqref{Eq_S2_asymptotics_7} and \eqref{Eq_S2_asymptotics_8} then finally gives the upper and lower bound \eqref{Eq_S2_estimate} of the Stirling numbers, namely
\begin{equation*}
\frac{c_F}{e^{\frac{1}{12}}}\frac{p!}{q!}g^S(p,q)e^{-qh^S(\frac{p}{q})}\leq S_2(p,q)\leq\frac{C_F}{2-e^{\frac{1}{12}}}\frac{p!}{q!}g^S(p,q)e^{-qh^S(\frac{p}{q})}.
\end{equation*}
ii)\;\;Since there is $\lim_{k\to\infty}\frac{m_k}{n_k}=:c>\frac{2}{3}$, we are able to choose $\widetilde{c}\in(\frac{2}{3},c)$ arbitrary but fixed. Corresponding to this $\widetilde{c}$, there exists some $\varepsilon>0$, and some $N\in\mathbb{N}$, such that
\begin{equation*}
\frac{\frac{1}{2}+\varepsilon}{\frac{1}{2}-\varepsilon}<\frac{3\widetilde{c}}{2},\qquad\text{and}\qquad\frac{m_k}{n_k}\geq\widetilde{c},\quad\text{for every }k\geq N.
\end{equation*}
Choose now $\mathcal{N}_0$ as the open ball of radius $\varepsilon$, centered at $(\frac{1}{2},\frac{1}{2})$. Then for every $(r,\alpha)\in\mathcal{N}_0$ and every $k\geq N$, there is
\begin{equation*}
\frac{\alpha}{r}\leq\frac{\frac{1}{2}+\varepsilon}{\frac{1}{2}-\varepsilon}<\frac{3\widetilde{c}}{2}\leq\frac{3m_k}{2n_k}.
\end{equation*}
Rearranging this inequality means that
\begin{equation}\label{Eq_S2_approximation_1}
\alpha n_k<\frac{3rm_k}{2},\qquad (r,\alpha)\in\mathcal{N}_0,\,k\geq N.
\end{equation}
Then the inequality \eqref{Eq_S2_asymptotics_3} with the values $p=3rm_k$ and $q=\alpha n_k$, gives
\begin{equation*}
\Big|\frac{F(3rm_k,\alpha n_k)}{S_2(3rm_k,\alpha n_k)}-1\Big|\leq E(3rm_k)\leq E\Big(3\Big(\frac{1}{2}-\varepsilon\Big)m_k\Big),\qquad(r,\alpha)\in\mathcal{N}_0\cap\Lambda_{m_k,n_k},\,k\geq N,
\end{equation*}
where in the second inequality we used $r\geq\frac{1}{2}-\varepsilon$ as well as the monotonicity of $E$. Note that in \cite[Theorem 1.2.1]{CHHarxiv}, the function $E$ is not assumed to be monotone decreasing, and only defined for integers $p\geq 1$. However, it is obvious that one can choose a function $E$ which satisfies the inequality \eqref{Eq_S2_asymptotics_3}, still converges as $\lim_{x\to\infty}E(x)=0$, is monotone decreasing as well as defined on the whole half line $(0,\infty)$, for example $\widetilde{E}(x):=\sup_{p\in\mathbb{N},p\geq x}E(p)$. Since the right hand side is independent of $(r,\alpha)$ this inequality also holds for the supremum
\begin{equation*}
\sup\limits_{(r,\alpha)\in\mathcal{N}_0\cap\Lambda_{m_k,n_k}}\Big|\frac{F(3rm_k,\alpha n_k)}{S_2(3rm_k,\alpha n_k)}-1\Big|\leq E\Big(3\Big(\frac{1}{2}-\varepsilon\Big)m_k\Big),\qquad k\geq N.
\end{equation*}
Finally, we can perform the limit $k\to\infty$ of this inequality. Using $\lim_{p\to\infty}E(p)=0$, one sees that the right hand side of this inequality vanishes, and we obtain the asymptotics
\begin{equation}\label{Eq_S2_asymptotics_2}
S_2(3rm_k,\alpha n_k)\sim F(3rm_k,\alpha n_k),\qquad\text{as }k\to\infty,\text{ uniformly in }(r,\alpha)\in\mathcal{N}_0\cap\Lambda_{m_k,n_k}.
\end{equation}
Since the Stirling approximation \Cref{Eq_S2_asymptotics_4} shows that
\begin{equation}\label{Eq_Factorial_asymptotics}
n!\sim\sqrt{2\pi n}\,\Big(\frac{n}{e}\Big)^n,\qquad\text{as }n\to\infty,
\end{equation}
we can further simplify the asymptotics \eqref{Eq_S2_asymptotics_2} to
\begin{align*}
S_2(3rm_k,\alpha n_k)&\sim F(3rm_k,\alpha n_k) \\
&=\frac{(3rm_k)!\sqrt{2\pi(3rm_k-2\alpha n_k)}}{(\alpha n_k)!(3rm_k-2\alpha n_k)!}\Big(\frac{3rm_k-2\alpha n_k}{e}\Big)^{3rm_k-2\alpha n_k} \\
&\hspace{5cm}\times g^S(3rm_k,\alpha n_k)e^{-\alpha n_kh^S(\frac{3rm_k}{\alpha n_k})} \\
&\sim\frac{(3rm_k)!}{(\alpha n_k)!}g^S(3rm_k,\alpha n_k)e^{-\alpha n_kh^S(\frac{3rm_k}{\alpha n_k})}. \qedhere
\end{align*}
\end{proof}

Let us now consider some crucial bound on the function $g^S$. The upcoming inequality \eqref{Eq_gS_estimate} is by no means obvious since $g^S$ has a singularity for $q\to(\frac{p}{2})^-$. However, since we only consider $(r,\alpha)\in\cW\cap\Lambda_{m,n}$ on the grid, it happens that there exists a uniform upper bound.

\begin{lemma}\label{lem_gS_estimate}
For every $p\in\mathbb{N}$, $q\in\mathbb{N}_0$ with $1\leq q<\frac{p}{2}$, there holds the estimate
\begin{equation}\label{Eq_gS_upper_lower_bound}
\frac{1}{\sqrt{2\pi p}}\leq g^S(p,q)\leq\frac{\sqrt{3}}{\sqrt{2\pi(p-2q)}}
\end{equation}
Furthermore, for every $m,n\in\mathbb{N}$ there holds the estimate
\begin{equation}\label{Eq_gS_estimate}
g^S(3rm,\alpha n)\leq 1,\qquad(r,\alpha)\in\widehat{\cW}\cap\Lambda_{m,n}.
\end{equation}
\end{lemma}

\begin{proof}
It can be checked easily from the definition of $Q$ in \eqref{Eq_Q}, that $Q^{-1}$ and $Q'$ admit the upper and lower bounds
\begin{equation*}
\xi-2\leq Q^{-1}(\xi)\leq\xi,\quad\xi\in(2,\infty),\qquad\text{and}\qquad\frac{1}{3}\leq Q'(x)\leq 1,\quad x\in(0,\infty).
\end{equation*}
From these two inequalities it then follows that
\begin{equation*}
\frac{\xi-2}{3}\leq Q^{-1}(\xi)Q'(Q^{-1}(\xi))\leq\xi,\qquad\xi\in(2,\infty),
\end{equation*}
and we immediately conclude that the function $g^S$ is for every $1\leq q<\frac{p}{2}$ bounded by
\begin{equation*}
\frac{1}{\sqrt{2\pi p}}\leq g^S(p,q)=\frac{1}{\sqrt{2\pi qQ^{-1}(\frac{p}{q})Q'(Q^{-1}(\frac{p}{q}))}}\leq\frac{\sqrt{3}}{\sqrt{2\pi(p-2q)}}.
\end{equation*}

For the second inequality \eqref{Eq_gS_estimate}, let $p=3rm$ and $q=\alpha n$. We only have to consider $1\leq q<\frac{p}{2}$, since all other cases are trivial by the definition of $g^S$ in \eqref{Eq_gS}. Since $p=3rm$ and $q=\alpha n$ are integers due to $(r,\alpha)\in\Lambda_{m,n}$, \Cref{Eq_gS_upper_lower_bound} holds with $p-2q\geq 1$. So we conclude the uniform upper bound
\begin{equation*}
g^S(3rm,\alpha n)\leq\frac{\sqrt{3}}{\sqrt{2\pi(3rm-2\alpha n)}}\leq\frac{\sqrt{3}}{\sqrt{2\pi}}\leq 1. \qedhere
\end{equation*}
\end{proof}

Next, we want to derive an estimate as well as an asymptotic formula for the binomial coefficient.
We will do this using the functions
\begin{align}
\mathdf{g^B(p,q)}:=&\frac{\sqrt{q}}{\sqrt{2\pi\max\{p,1\}\max\{q-p,1\}}},\qquad p\in\mathbb{N}_0,\,q\in\mathbb{N},\,p\leq q, \label{Eq_gB} \\
\mathdf{h^B(\xi)}:=&\begin{cases} \xi\ln(\xi)+(1-\xi)\ln(1-\xi), & \xi\in(0,1), \\ 0, & \xi=0\text{ or }\xi=1. \end{cases} \label{Eq_hB}
\end{align}
Note that the function $h^B$ is continuous in $\xi=0$ and in $\xi=1$.

\begin{lemma}
i)\;\;There exist constants $0<c_B\leq C_B$ such that for every $p\in\mathbb{N}_0$, $q\in\mathbb{N}$ with $p\leq q$, the binomial coefficient is bounded from below and from above by
\index{$c_B$, $C_B$}
\begin{equation}\label{Eq_Binom_estimate}
c_B\,g^B(p,q)e^{-qh^B(\frac{p}{q})}\leq\binom{q}{p}\leq C_B\,g^B(p,q)e^{-qh^B(\frac{p}{q})}.
\end{equation}
ii)\;\;Let $(p_k)_k,(q_k)_k\subseteq\mathbb{N}$, with $\lim_{k\to\infty}p_k=\lim_{k\to\infty}q_k=\infty$ and $\lim_{k\to\infty}\frac{p_k}{q_k}\in(0,1)$. Then the binomial coefficient admits the asymptotics
\begin{equation}\label{Eq_Binom_asymptotics}
\binom{q_k}{p_k}\sim g^B(p_k,q_k)e^{-q_kh^B(\frac{p_k}{q_k})},\qquad\text{as }k\to\infty.
\end{equation}
\end{lemma}

\begin{proof}
i)\;\;For $p=0$ or $p=q$, there clearly holds the equality
\begin{equation*}
\binom{q}{p}=1=\frac{\sqrt{q}}{\sqrt{\max\{p,1\}\max\{q-p,1\}}}=\sqrt{2\pi}\,g^B(p,q)e^{-qh^B(\frac{p}{q})}.
\end{equation*}
So we only have to consider $0<p<q$. Using the estimate \eqref{Eq_Factorial_estimate} three times, with the values $n=q$, $n=p$ and $n=q-p$ respectively, gives the following upper and lower bounds of the binomial coefficient
\begin{equation*}
\frac{(2-e^{\frac{1}{12}})\sqrt{q}\,(\frac{p}{q})^{-p}(1-\frac{p}{q})^{-(q-p)}}{e^{\frac{1}{6}}\sqrt{2\pi p(q-p)}}\leq\binom{q}{p}=\frac{q!}{p!(q-p)!}\leq\frac{e^{\frac{1}{12}}\sqrt{q}\,(\frac{p}{q})^{-p}(1-\frac{p}{q})^{-(q-p)}}{(2-e^{\frac{1}{12}})^2\sqrt{2\pi p(q-p)}}.
\end{equation*}
With the definitions \eqref{Eq_gB} and \eqref{Eq_hB} of $g^B$ and $h^B$, we can then rewrite this inequality as
\begin{equation*}
\frac{2-e^{\frac{1}{12}}}{e^{\frac{1}{6}}}g^B(p,q)e^{-qh^B(\frac{p}{q})}\leq\binom{q}{p}\leq\frac{e^{\frac{1}{12}}}{(2-e^{\frac{1}{12}})^2}g^B(p,q)e^{-qh^B(\frac{p}{q})}.
\end{equation*}
ii)\;\;Using the Stirling asymptotics \eqref{Eq_Factorial_asymptotics} three times, with the values $n_k=q_k$, $n_k=p_k$ and $n_k=q_k-p_k$ respectively, gives the asymptotics of the binomial coefficient
\begin{equation*}
\binom{q_k}{p_k}=\frac{q_k!}{p_k!(q_k-p_k)!}\sim\frac{\sqrt{q_k}\,(\frac{p_k}{q_k})^{-p_k}(\frac{q_k-p_k}{q_k})^{-(q_k-p_k)}}{\sqrt{2\pi p_k(q_k-p_k)}}=g^B(p_k,q_k)e^{-q_kh^B(\frac{p_k}{q_k})}. \qedhere
\end{equation*}
\end{proof}

\subsection{Asymptotics of the full summand $\mathcal{S}_{m,n}(r,\alpha)$}\label{sec_Asymptotics_of_Smn}

The results of \Cref{sec_Ingredients_of_Smn} will now be used to estimate the full summand $\mathcal{S}_{m,n}(r,\alpha)$ in Eq.\eqref{Eq_Smn_appendix}. As we shall see in the proof (see \Cref{sec_Subsection3}), the limit in Eq.\eqref{Eq_Application_3XORSAT_appendix} is the same if $r=1$ is excluded from the sum. so our preliminary lemmas do not need to address these cases. To this end, we define
\begin{equation*}
\cW\coloneqq\big\{(r,\alpha)\in\hW\;\big|\;r\neq 1\big\}=[\mathsmaller{\frac{1}{3}},1)\times[0,1],
\end{equation*}
and use it instead of $\widehat{\cW}$. This avoids a singularity of the upcoming function $g_{m,n}(r,\alpha)$ at $(r,\alpha)=(1,1)$. Using the functions $g^S$ and $g^B$ from \eqref{Eq_gS} and \eqref{Eq_gB}, we introduce for every $m,n\in\mathbb{N}$ with $2<\frac{3m}{n}<3$, the function $g_{m,n}\colon\mathcal{W}\cap\Lambda_{m,n}\rightarrow\mathbb{R}$ by
\begin{equation}\label{Eq_gmn}
\mathdf{g_{m,n}(r,\alpha)}:=n\frac{g^S(3rm,\alpha n)g^S(3(1-r)m,(1-\alpha)n)g^B(\alpha n,n)g^B(\frac{3}{2}(1-r)m,m)}{g^S(3m,n)g^B(3rm,3m)}.
\end{equation}
Moreover, using $h^S$ and $h^B$ from \eqref{Eq_hS} and \eqref{Eq_hB}, we define for every $2<y< 3$ the function $h_y:\cW\rightarrow\mathbb{R}$ by
\begin{equation}\label{Eq_hy}
\mathdf{h_y(r,\alpha)}:=\begin{cases} \alpha h^S\Big(\frac{ry}{\alpha}\Big)+(1-\alpha)h^S\Big(\frac{(1-r)y}{1-\alpha}\Big)-h^S(y)+h^B(\alpha)+\ln(2) \\
\qquad-y\Big(h^B(r)-\frac{1}{3}h^B\Big(\frac{3(1-r)}{2}\Big)+\frac{\ln(2)}{3}+\frac{\ln(3)(1-r)}{2}\Big), & (r,\alpha)\in\mathcal{T}_y, \\ 1, & (r,\alpha)\notin\mathcal{T}_y. \end{cases}
\end{equation}

\begin{figure}[H]
\begin{minipage}{.4\textwidth}
\centering
\begin{tikzpicture}
\fill[black!25] (4,4)--(3.08,4)--(1.33,1.73)--(1.33,0.45);
\draw[->] (0,0)--(4.5,0) node[anchor=south] {$r$};
\draw[->] (0,0)--(0,4.4) node[anchor=west] {$\alpha$};
\draw[ultra thick] (4,4)--(1.33,4)--(1.33,0)--(4,0);
\draw[dashed, ultra thick] (4,0)--(4,4);
\draw (1.33,0) node[anchor=north] {$\frac{1}{3}$};
\draw (4,0) node[anchor=north] {$1$};
\draw (-0.1,0) node[anchor=east] {$0$};
\draw (0.1,4)--(-0.1,4) node[anchor=east] {$1$};
\draw[dashed] (0,0)--(1.33,1.73);
\draw[ultra thick] (1.33,1.73)--(3.08,4) node[anchor=south] {\tiny{$\frac{2}{y}$}};
\draw[dashed] (0.92,0)--(1.33,0.45);
\draw[ultra thick] (1.33,0.45)--(4,4);
\draw (0.7,0) node[anchor=north] {\tiny{$1$-$\frac{2}{y}$}};
\fill[black] (2,2) circle (0.07cm);
\draw (2.2,2) node[anchor=south] {\tiny{$(\frac{1}{2},\frac{1}{2})$}};
\draw (3.1,1) node[anchor=south] {\large{$\cW$}};
\draw (2.8,3.4) node[anchor=north] {\large{$\cT_y$}};
\draw[thick,fill=white] (4,4) circle (0.06cm);
\draw[thick,fill=white] (4,0) circle (0.06cm);
\end{tikzpicture}
\end{minipage}
\begin{minipage}{.5\textwidth}
\centering
\begin{align*}
\mathdf{\cT_y} \coloneqq \Big\{(r,\alpha) \in \cW\;\Big|\; 1 - \frac{(1-r)y}{2} \leq \alpha\leq \frac{ry}{2} \Big\}
\end{align*}
\end{minipage}
\caption{Illustration of the region $\cT_y$ where for $y=\frac{3m}{n}$ the function $g_{m,n}(r,\alpha)$ is nonzero. Outside $\cT_y$, the function $g_{m,n}(r,\alpha)$ vanishes and $h_{3m/n}(r,\alpha)$ is set identically one (actually, any positive number would do). The interior of $\cT_y$ is also the region where $g_{m,n}(r,\alpha)$ coincides with $g_{3m/n}(r,\alpha)$ in Eq.\eqref{Eq_g}, and also where the product $g_{m,n}(r,\alpha)e^{-nh_{3m/n}(r,\alpha)}$ asymptotically behaves as the summand $\mathcal{S}_{m,n}(r,\alpha)$, see Eq.\eqref{Eq_Smn_asymptotics}.  On the boundary of $\cT_y$, the value of $g_{m,n}(r,\alpha)e^{-nh_{3m/n}(r,\alpha)}$ is set such that it is an upper bound of the summand $\mathcal{S}_{m,n}(r,\alpha)$, see Eq.\eqref{Eq_Smn_estimate}.}
\label{fig_Ty}
\end{figure}
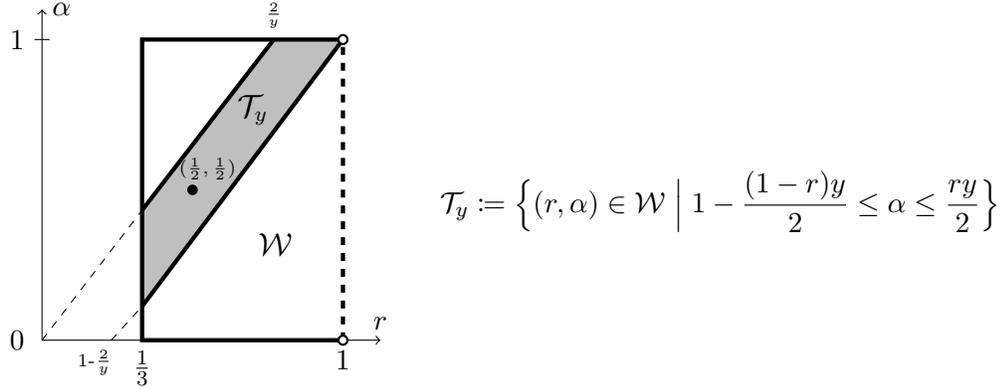

\begin{theorem}
i)\;\;Let $m,n\in\mathbb{N}$ with $2<\frac{3m}{n}<3$. Then, with the constants from \eqref{Eq_S2_estimate} and \eqref{Eq_Binom_estimate}, the summand \eqref{Eq_Smn} satisfies the estimate
\begin{equation}\label{Eq_Smn_estimate}
\mathcal{S}_{m,n}(r,\alpha)\leq\frac{C_S^2C_B^2}{c_Sc_B}g_{m,n}(r,\alpha)e^{-nh_{3m/n}(r,\alpha)},\qquad(r,\alpha)\in\cW\cap\Lambda_{m,n}
\end{equation}
ii)\;\;Let $(m_k)_k,(n_k)_k\subseteq\mathbb{N}$, with $\lim_{k\to\infty}m_k=\lim_{k\to\infty}n_k=\infty$ and $\lim_{k\to\infty}\frac{m_k}{n_k}\in(\frac{2}{3},1)$. Then there exists an open neighborhood $\mathcal{N}_0$ of $(\frac{1}{2},\frac{1}{2})$, such that
\begin{equation}\label{Eq_Smn_asymptotics}
\mathcal{S}_{m_k,n_k}(r,\alpha)\sim g_{m_k,n_k}(r,\alpha)e^{-n_kh_{3m_k/n_k}(r,\alpha)},
\end{equation}
as $k\to\infty$, uniformly for $(r,\alpha)\in\mathcal{N}_0\cap\Lambda_{m_k,n_k}$, in the sense of \eqref{Eq_Asymptotics_uniform}.
\end{theorem}

\begin{proof}
i)\;\; First note that by the definition of the Stirling numbers in Eq.\eqref{Eq_Stirling_number}, there is
\begin{align*}
S_2(3rm,\alpha n)&=0,\qquad\text{for }\alpha>\frac{3rm}{2n}, \\
S_2(3(1-r)m,(1-\alpha)n)&=0,\qquad\text{for }\alpha<1-\frac{3(1-r)m}{2n}.
\end{align*}
This in particular means that $S_{m,n}(r,\alpha)=0$ for every $(r,\alpha)\notin\mathcal{T}_{3m/n}$, and hence the inequality \eqref{Eq_Smn_estimate} is trivially satisfied.
\medskip

For $(r,\alpha)\in\mathcal{T}_{3m/n}$ on the other hand, the inequality \eqref{Eq_S2_estimate} of the Stirling numbers allows us to estimate \eqref{Eq_Smn_appendix} by
\begin{align*}
\mathcal{S}_{m,n}(r,\alpha)&\leq\frac{C_S^2}{c_S}n2^{m-n}3^{\frac{3}{2}(1-r)m}\binom{m}{\frac{3}{2}(1-r)m}\binom{3m}{3rm}^{-1}\binom{n}{\alpha n} \\
&\quad\times\frac{g^S(3rm,\alpha n)g^S(3(1-r)m,(1-\alpha)n)}{g^S(3m,n)}e^{-\alpha nh^S(\frac{3rm}{\alpha n})}e^{-(1-\alpha)nh^S(\frac{3(1-r)m}{(1-\alpha)n})}e^{nh^S(\frac{3m}{n})}.
\end{align*}
Applying also the estimate \eqref{Eq_Binom_estimate} of the binomial coefficients, then gives the stated inequality
\begin{align*}
\mathcal{S}_{m,n}(r,\alpha)&\leq\frac{C_S^2C_B^2}{c_Sc_B}n2^{m-n}3^{\frac{3}{2}(1-r)m}\frac{g^B(\alpha n,n)g^B(\frac{3}{2}(1-r)m,m)}{g^B(3rm,3m)} \\
&\qquad\times\frac{g^S(3rm,\alpha n)g^S(3(1-r)m,(1-\alpha)n)}{g^S(3m,n)} \\
&\qquad\times e^{-\alpha nh^S(\frac{3rm}{\alpha n})}e^{-(1-\alpha)nh^S(\frac{3(1-r)m}{(1-\alpha)n})}e^{nh^S(\frac{3m}{n})}e^{-nh^B(\alpha)}e^{3mh^B(r)}e^{-mh^B(\frac{3}{2}(1-r))} \\
&=\frac{C_S^2C_B^2}{c_Sc_B}g_{m,n}(r,\alpha)e^{-nh_{3m/n}(r,\alpha)}.
\end{align*}
ii)\;\;In view of the approximation \eqref{Eq_S2_approximation}, there exists a neighborhood $\mathcal{N}_0$ of $(\frac{1}{2},\frac{1}{2})$, such that the right hand side of Eq.\eqref{Eq_Smn} asymptotically behaves as
\begin{align*}
\mathcal{S}_{m_k,n_k}(r,\alpha)&\sim n_k2^{m_k-n_k}3^{\frac{3}{2}(1-r)m_k}\binom{m_k}{\frac{3}{2}(1-r)m_k}\binom{3m_k}{3rm_k}^{-1}\binom{n_k}{\alpha n_k} \\
&\qquad\times\frac{g^S(3rm_k,\alpha n_k)g^S(3(1-r)m_k,(1-\alpha)n_k)}{g^S(3m_k,n_k)} \\
&\qquad\times e^{-\alpha n_kh^S(\frac{3rm_k}{\alpha n_k})}e^{-(1-\alpha)n_kh^S(\frac{3(1-r)m_k}{\alpha n_k})}e^{n_kh^S(\frac{3m_k}{n_k})},
\end{align*}
as $k\to\infty$, uniformly for $(r,\alpha)\in\mathcal{N}_0\cap\Lambda_{m_k,m_k}$. The asymptotic \eqref{Eq_Binom_asymptotics} of the binomial coefficient allows us to further expand the right-hand side to
\begin{align*}
\mathcal{S}_{m_k,n_k}(r,\alpha)&\sim n_k2^{m_k-n_k}3^{\frac{3}{2}(1-r)m_k}\frac{g^B(\frac{3}{2}(1-r)m_k,m_k)g^B(\alpha n_k,n_k)}{g^B(3rm_k,m_k)} \\
&\qquad\times\frac{g^S(3rm_k,\alpha n_k)g^S(3(1-r)m_k,(1-\alpha)n_k)}{g^B(3rm_k,m_k)} \\
&\qquad\times e^{-\alpha n_kh^S(\frac{3rm_k}{\alpha n_k})}e^{-(1-\alpha)n_kh^S(\frac{3(1-r)m_k}{\alpha n_k})}e^{3m_kh^B(r)}e^{-m_kh^B(\frac{3}{2}(1-r))}e^{-n_kh^B(\alpha)} \\
&=g_{m_k,n_k}(r,\alpha)e^{-n_kh_{3m_k/n_k}(r,\alpha)},
\end{align*}
as $k\to\infty$, uniformly for $(r,\alpha)\in\mathcal{N}_0\cap\Lambda_{m_k,n_k}$.
\end{proof}

\subsection{Positivity of the function $h_y$}\label{sec_Positivity_hy}

The aim of this subsection is to show that for all $y\in(2,3)$, the global minimum value of the function $h_y$ in \eqref{Eq_hy} is attained at the point $(r,\alpha)=(\frac{1}{2},\frac{1}{2})$. This is the statement of \cite[Lemma 3.5]{dBM-FOCS}, and we present a proof here. Our emphasis is on filling in details left to the reader in their proof. Such omission is primarily checking positivity on the boundary and this seems tricky, for example, $h_y$ is not even continuous at $(r,\alpha)=(1,1)$.

\begin{theorem}\label{thm_hy_Minimum}
For any $y\in(2,3)$, there is
\begin{equation}\label{Eq_hy_infimum}
\inf\limits_{(r,\alpha)\in\cW}h_y(r,\alpha)=0,
\end{equation}
and the infimum is attained only at the point $(r,\alpha)=(\frac{1}{2},\frac{1}{2})$.
\end{theorem}

The proof will be divided into two parts; each is a subsection, and the final Subsection~\ref{sec:conclude-lemma-hy-W}, draws on these to complete the proof.

\subsubsection{Properties of $h_y(r, \alpha)$ as a function of $\alpha$}

For $y\in(2,3)$ and fixed $r\in{[}\frac{1}{3},1)$, we will first investigate the mapping $\alpha\mapsto h_y(r,\alpha)$, defined for every $\alpha$ such that $(r,\alpha)\in\mathcal{T}_y$. Note, as it can be seen in \Cref{fig_Ty}, the condition $(r,\alpha)\in\mathcal{T}_y$ is, for fixed $r$, equivalent to
\begin{equation}\label{Eq_alpha_interval}
1-\frac{(1-r)y}{2}\leq\alpha\leq\min\Big\{\frac{ry}{2},1\Big\}.
\end{equation}

\begin{lemma}\label{lem_hy_convex}
For fixed $y\in(2,3)$ and $r\in{[}\frac{1}{3},1)$, the mapping $\alpha\mapsto h_y(r,\alpha)$ is strictly convex.
\end{lemma}

\begin{proof} 
Using the derivatives
\begin{equation*}
\frac{d}{d\xi}h^S(\xi)=\ln(Q^{-1}(\xi))\qquad\text{and}\qquad\frac{d}{d\xi}h^B(\xi)=\ln\Big(\frac{\xi}{1-\xi}\Big),
\end{equation*}
of \eqref{Eq_hS} and \eqref{Eq_hB}, we can calculate the $\alpha$-derivative of the function $h_y$ in \eqref{Eq_hy}
\begin{align}
\mathdf{\frac{\partial}{\partial\alpha}h_y(r,\alpha)}&=h^S\Big(\frac{ry}{\alpha}\Big)-\frac{ry}{\alpha}{h^S}'\Big(\frac{ry}{\alpha}\Big)-h^S\Big(\frac{(1-r)y}{1-\alpha}\Big)+\frac{(1-r)y}{1-\alpha}{h^S}'\Big(\frac{(1-r)y}{1-\alpha}\Big)+{h^B}'(\alpha) \notag \\
&=-\ln\Big(A\Big(\frac{ry}{\alpha}\Big)\Big)+\ln\Big(A\Big(\frac{(1-r)y}{1-\alpha}\Big)\Big)+\ln\Big(\frac{\alpha}{1-\alpha}\Big). \label{eq:derivAlpha}
\end{align}
where for convenience we used once more the function $A(\xi)$ in \eqref{Eq_A}. Differentiating this expression once more with respect to $\alpha$, then gives
\begin{equation*}
\frac{\partial^2}{\partial\alpha^2}h_y(r,\alpha)=\frac{ry A'(\frac{ry}{\alpha})}{\alpha^2A(\frac{ry}{\alpha})}+\frac{(1-r)yA'(\frac{(1-r)y}{1-\alpha})}{(1-\alpha)^2A(\frac{(1-r)y}{1-\alpha})}+\frac{1}{\alpha(1-\alpha)}.
\end{equation*}
Since it is clear that $A(\xi)>0$ and $A'(\xi)>0$ for every $\xi>2$, every single term of this expression is positive. Consequently, $\frac{\partial^2}{\partial\alpha^2}h_y(r,\alpha)>0$ is positive, which means that for fixed $r\in{[}1/3,1)$, the mapping $\alpha\mapsto h_y(r,\alpha)$ is strictly convex.
\end{proof}

\begin{lemma}\label{lem_alphay}
Let $y\in(2,3)$ and $r\in{[}\frac{1}{3},1)$. Then there exists a unique minimizer $\alpha_y(r)$ of
\begin{equation}\label{Eq_hy_Minimization}
\inf\limits_{\alpha\in[0,1]: (r,\alpha)\in\mathcal{T}_y}h_y(r,\alpha).
\end{equation}
This minimizer satisfies $(r,\alpha_y(r))\in\Int(\mathcal{T}_y)$, is the unique solution of the equation
\begin{equation}\label{Eq_alphay_definition}
\frac{\alpha A(\frac{(1-r)y}{1-\alpha})}{(1-\alpha)A(\frac{ry}{\alpha})}=1,
\end{equation}
and is located at
\begin{equation}\label{Eq_alphay_location}
\mathdf{\alpha_y(r)} \begin{cases} >r, & \text{if }\frac{1}{3}\leq r<\frac{1}{2}, \\ =r, & \text{if }r=\frac{1}{2}, \\ <r, & \text{if }\frac{1}{2}<r<1. \end{cases}
\end{equation}
\end{lemma}

\begin{figure}[H]
\begin{center}
\begin{tikzpicture}
\fill[black!25] (4,4)--(3.08,4)--(1.33,1.73)--(1.33,0.45);
\draw[->] (0,0)--(4.5,0) node[anchor=south] {$r$};
\draw[->] (0,0)--(0,4.4) node[anchor=west] {$\alpha$};
\draw[thick] (4,4)--(1.33,4)--(1.33,0)--(4,0);
\draw[dashed, thick] (4,0)--(4,4);
\draw (1.33,0) node[anchor=north] {$\frac{1}{3}$};
\draw (4,0) node[anchor=north] {$1$};
\draw (-0.1,0) node[anchor=east] {$0$};
\draw (0.1,4)--(-0.1,4) node[anchor=east] {$1$};
\draw[dashed] (0,0)--(1.33,1.73);
\draw[thick] (1.33,1.73)--(3.08,4) node[anchor=south] {\tiny{$\frac{2}{y}$}};
\draw[dashed] (0.92,0)--(1.33,0.45);
\draw[thick] (1.33,0.45)--(4,4);
\draw (0.7,0) node[anchor=north] {\tiny{$1$-$\frac{2}{y}$}};
\fill[black] (1.98,2) circle (0.07cm);
\draw (1.75,2) node[anchor=south] {\tiny{$(\frac{1}{2},\frac{1}{2})$}};
\draw (3.1,1) node[anchor=south] {\large{$\cW$}};
\draw (3,3.1) node[anchor=south] {\scriptsize{$\alpha_y(r)$}};
\draw[ultra thick] plot[domain=0.34:0.68] (4*\x-0.02,{(0.067-0.133*\x+\x)*4});
\draw[ultra thick] plot[domain=0.68:0.81] (4*\x-0.02,{(0.4545*(\x-0.81)^2-0.031+\x)*4});
\draw[ultra thick] plot[domain=0.81:1] (4*\x-0.02,{(0.859*(\x-0.81)^2-0.031+\x)*4});
\draw[thick,fill=white] (4,4) circle (0.06cm);
\draw[thick,fill=white] (4,0) circle (0.06cm);
\end{tikzpicture}
\end{center}
\caption{The heavy curve through the middle of the trapezoid denote the points $(r, \alpha_y(r) )$, given by the implicit equation \eqref{Eq_alphay_definition}. These are the unique minimizers of the function $\alpha\mapsto h_y(r,\alpha)$. This is also the curve along which we will minimize the function $h_y$ in the second step of the proof, in order to find the global minimum.}
\end{figure}
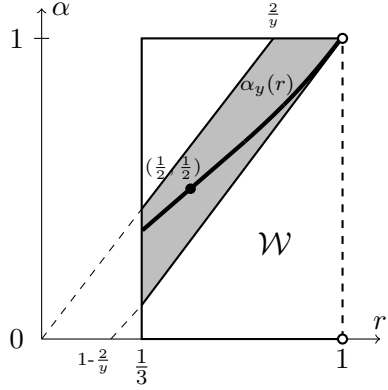

\begin{proof}
From \Cref{lem_hy_convex} we know that the mapping $\alpha\mapsto h_y(r,\alpha)$ is strictly convex. Hence there exists unique minimizer $\alpha_y(r)$ of the infimum \eqref{Eq_hy_Minimization}, possibly located at the endpoints of the allowed interval \eqref{Eq_alpha_interval} of $\alpha$. To understand the minimizer $\alpha_y(r)$ when it is in the
interior, we set the derivative formula \eqref{eq:derivAlpha} to zero. The equation this gives for $\alpha$, is
\begin{equation*}
\mathdf{L_{r,y}(\alpha)}:=\frac{\alpha A\big(\frac{(1-r)y}{1-\alpha}\big)}{(1-\alpha)A\big(\frac{ry}{\alpha}\big)}=1,\qquad 1-\frac{(1-r)y}{2}<\alpha<\min\Big\{\frac{ry}{2},1\Big\}.
\end{equation*}
Since the function $A$ in \eqref{Eq_A} is monotone increasing, the complete left hand side $L_{r,y}(\alpha)$ is monotone increasing as a function in $\alpha$. Moreover, at the boundary points $\alpha=1-\frac{(1-r)y}{2}$ and $\alpha=\min\{\frac{ry}{2},1\}$, the function $L_{r,y}$ admits the values
\begin{align*}
\lim\limits_{\alpha\rightarrow 1-\frac{(1-r)y}{2}}L_{r,y}(\alpha)&=\frac{2-(1-r)y}{(1-r)yA(\frac{2ry}{2-(1-r)y})}\lim\limits_{\xi\rightarrow 2^+}A(\xi)=0, \\
\lim\limits_{\alpha\rightarrow\min\{\frac{ry}{2},1\}}L_{r,y}(\alpha)&=\left.\begin{cases} \frac{ry}{2-ry}A(\frac{2(1-r)y}{2-ry})\lim\limits_{\xi\rightarrow 2^+}\frac{1}{A(\xi)}, & \text{if }\frac{ry}{2}<1, \\ \frac{1}{(1-r)yA(ry)}\lim\limits_{\xi\rightarrow\infty}\xi A(\xi), & \text{if }\frac{ry}{2}\geq 1, \end{cases}\right\} =\infty.
\end{align*}
These boundary values together with the monotonicity of $L_{r,y}$ show, that the unique minimizer $\alpha_y(r)$ exists in the interior of the interval $1-\frac{(1-r)y}{2}<\alpha_y(r)<\min\{\frac{ry}{2},1\}$. In other words, this means $(r,\alpha_y(r))\in\Int(\mathcal{T}_y)$. In order to find the location \eqref{Eq_alphay_location} of $\alpha_y(r)$ relative to $r$, we evaluate $L_{r,y}(\alpha)$ at $\alpha=r$, which yields
\begin{equation*}
L_{r,y}(r)=\frac{r}{1-r}\begin{cases} < 1, & \text{if }\frac{1}{3}\leq r<\frac{1}{2}, \\ =1, &  \text{if }r=\frac{1}{2}, \\ >1, &  \text{if }\frac{1}{2}<r<1. \end{cases}
\end{equation*}
Since $\alpha\mapsto L_{r,y}(\alpha)$ is monotone increasing, and $\alpha_y(r)$ is the unique root to $L_{r,y}(\alpha)=1$, this then implies the location \eqref{Eq_alphay_location} of $\alpha_y(r)$ relative to the line $\alpha=r$.
\end{proof}

\subsubsection{The minimal values $h_y(r,\alpha_y(r))$ as a function of $r$.}

Now that we have calculated the function $\alpha_y(r)$ in \Cref{lem_alphay}, we know that the minimal value of $h_y(r,\alpha)$ is located somewhere on the line $(r,\alpha_y(r))$, $r\in{[}\frac{1}{3},1)$. Hence we have reduced the problem \eqref{Eq_hy_infimum} to the one dimensional infimum of the function $h_y(r,\alpha_y(r))$, as a function of $r$. Let us start with the behavior of this function on the left and on the right endpoints $r=\frac{1}{3}$ and $r\to 1$.

\begin{lemma}\label{lem:hy-near-1}
Let $y\in(2,3)$. Then with the function $\alpha_y(r)$ from \Cref{lem_alphay}, there is \medskip

\begin{enumerate}
\item[i)] $h_y\big(\frac{1}{3},\alpha_y(\frac{1}{3})\big)\geq y\big(\frac{2\ln(3)}{3}-\ln(2)\big)>0$, \medskip

\item[ii)] $\lim\limits_{r\to 1}h_y(r,\alpha_y(r))=\ln(2)(1-\frac{y}{3})>0$.
\end{enumerate}
\end{lemma}

\begin{proof}
i)\;\;First, it is easy to check that the $h^B$ in \eqref{Eq_hB} is bounded from below by
\begin{equation*}
h^B(\xi)\geq-\ln(2),\qquad\xi\in[0,1].
\end{equation*}
Using this, we can estimate for every $\alpha\in[0,1]$ with $(\frac{1}{3},\alpha)\in\mathcal{T}_y$, the function $h_y(\frac{1}{3},\alpha)$ by
\begin{align*}
h_y\Big(\frac{1}{3},\alpha\Big)&=\alpha h^S\Big(\frac{y}{3\alpha}\Big)+(1-\alpha)h^S\Big(\frac{2y}{3(1-\alpha)}\Big)-h^S(y)+h^B(\alpha)+(1-y)\ln(2)+\frac{2y\ln(3)}{3}  \\
&\geq\alpha h^S\Big(\frac{y}{3\alpha}\Big)+(1-\alpha)h^S\Big(\frac{2y}{3(1-\alpha)}\Big)-h^S(y)+y\Big(\frac{2\ln(3)}{3}-\ln(2)\Big) \\
&\geq y\Big(\frac{2\ln(3)}{3}-\ln(2)\Big)>0,
\end{align*}
where in the last line we used convexity of $h^S(\xi)$, which follows immediately from the positivity of its second derivative \begin{equation*}
\frac{d^2}{d\xi^2}h^S(\xi)=\frac{1}{Q^{-1}(\xi)Q'(Q^{-1}(\xi))}>0.
\end{equation*}
Since this is true for every $\alpha$, it is in particular true for $\alpha=\alpha_y(\frac{1}{3})$. \medskip

ii)\;\; 
While the limit of $h_y(r,\alpha)$ as $(r,\alpha)\to(1,1)$ exists along various curves, this limit along different curves may produce different values. However, we will show that the limit along the curve $(r,\alpha_y(r))$ exists and gives the stated value. Plugging in the value $\alpha=\alpha_y(r)$ into the definition \eqref{Eq_hy} of the function $h_y$ gives
\begin{align*}
h_y(r,\alpha_y(r))&=\alpha_y(r)h^S\Big(\frac{ry}{\alpha_y(r)}\Big)+(1-\alpha_y(r))h^S\Big(\frac{(1-r)y}{1-\alpha_y(r)}\Big)-h^S(y)+h^B(\alpha_y(r))+\ln(2) \\
&\quad-y\Big(h^B(r)-\frac{1}{3}h^B\Big(\frac{3(1-r)}{2}\Big)+\frac{\ln(2)}{3}+\frac{\ln(3)(1-r)}{2}\Big).
\end{align*}
Performing the limit $r\to 1$, and using $\lim_{r\to 1}\alpha_y(r)=1$ as well as 
$h^B(0)=h^B(1)=0$, we get
\begin{equation}\label{Eq_Left_boundary_1}
\lim\limits_{r\to 1}h_y(r,\alpha_y(r))=\lim\limits_{r\to 1}(1-\alpha_y(r))h^S\Big(\frac{(1-r)y}{1-\alpha_y(r)}\Big)+\ln(2)\Big(1-\frac{y}{3}\Big).
\end{equation}
For the remaining limit on the right hand side, we note that by \eqref{Eq_alphay_location} and the shape of the trapezoid $\mathcal{T}_y$, the function $\alpha_y(r)$ is bounded by $1-\frac{(1-r)y}{2}\leq\alpha_y(r)\leq r$, for every $\frac{1}{2}\leq r\leq 1$. Consequently, the term
\begin{equation*}
\frac{2}{y}\leq\frac{1-r}{1-\alpha_y(r)}\leq 1.\qquad\frac{1}{2}\leq r<1,
\end{equation*}
is bounded, and so $h^S\big(\frac{(1-r)y}{1-\alpha_y(r)}\big)$ is bounded as well. This means the limit \eqref{Eq_Left_boundary_1} reduces to
\begin{equation*}
\lim\limits_{r\to 1}h_y(r,\alpha_y(r))=\ln(2)\Big(1-\frac{y}{3}\Big). \qedhere
\end{equation*}
\end{proof}

\begin{proposition}\label{prop:hy-on-alpha-y}
Let $y\in(2,3)$. Then with $\alpha_y(r)$ from \Cref{lem_alphay}, there is
\begin{equation*}
\inf\limits_{r\in[\frac{1}{3},1)}h_y(r,\alpha_y(r))=0,
\end{equation*}
and the infimum is attained only at the point $r=\frac{1}{2}$.
\end{proposition}

\begin{proof}
\Cref{lem:hy-near-1} i) shows that at the left endpoint $r=\frac{1}{3}$, there is
\begin{equation*}
h_y(\mathsmaller{\frac{1}{3}},\alpha_y(\mathsmaller{\frac{1}{3}}))>0.
\end{equation*}
\Cref{lem:hy-near-1} ii) on the other hand shows that also on the right endpoint $r\to 1$, there is
\begin{equation*}
\lim_{r\to 1}h_y(r,\alpha_y(r))>0.
\end{equation*}
So the infimum cannot be located on either of the two endpoints, thus it must be a local minimum attained for some $r\in(\frac{1}{3},1)$. However, Dubois-Mandler in \cite[Lemma 3.5]{dBM-FOCS} prove that $h_y$ has only one local minimum on the interior of $\cT_y$; located at $(\half,\half)$ and with
\begin{equation*}
h_y(\mathsmaller{\frac{1}{2},\frac{1}{2}})=0.
\end{equation*}
Their proof has considerable detail so we do not reproduce it. (Note the notational difference; they instead prove $f_c\leq 2(1-c)\ln(2)$, where $f_c(r,\alpha):=2(1-\frac{y}{3})\ln(2)-h_y(r,\alpha)$ and $c=\frac{y}{3}$). Since $\alpha_y(\half)=\half$ by \eqref{Eq_alphay_location}, this implies the infimum of $h_y(r,\alpha_y(r))$ is attained only at $r=\half$.
\end{proof}

\subsubsection{Proof of \Cref{thm_hy_Minimum}}\label{sec:conclude-lemma-hy-W}

Now that all the preparations are done, we are in the position to prove the initial theorem of this subsection.

\begin{proof}[Proof of \Cref{thm_hy_Minimum}]
The goal is to show $h_y(r,\alpha)>0$ for all $(r,\alpha)\in\cW\setminus\{(\frac{1}{2},\frac{1}{2})\}$. To start proving this, note that on the exterior of the trapezoid, the value of $h_y$ is given by
\begin{equation*}
h_y(r,\alpha)=1,\qquad(r,\alpha)\in\cW\setminus\cT_y.
\end{equation*}
Thus, the problem is reduced to finding the minimum of
$h_y$ in the trapezoid $\cT_y$. \medskip

\Cref{lem_alphay} establishes a function $\alpha_y(r)$, such that the infimum of $h_y(r,\alpha)$, $(r,\alpha)\in\mathcal{T}_y$, must be located somewhere on the curve $(r,\alpha_y(r))$. Thus the problem reduces to finding the infimum of $h_y(r,\alpha_y(r))$ for $r\in{[}\frac{1}{3},1)$. However, by \Cref{prop:hy-on-alpha-y}, this infimum is in fact attained only at $r=\frac{1}{2}$, with
\begin{equation*}
h_y(\mathsmaller{\frac{1}{2}},\alpha_y(\mathsmaller{\frac{1}{2}}))=h_y(\mathsmaller{\frac{1}{2},\frac{1}{2}})=0. \qedhere
\end{equation*}
\end{proof}

\subsubsection{Uniform positivity of $h_y$}

One weakness of \Cref{thm_hy_Minimum} is that it only considers fixed values of $y$ when taking the infimum over the function $h_y(r,\alpha)$. However, in assumption \eqref{Eq_hk_lower_bound}, we need some kind of uniformity of the infimum after removing a neighborhood $\cN$ of $(\frac{1}{2},\frac{1}{2})$. The upcoming corollary handles this.

\begin{corollary}\label{cor_hy_liminf}
Let $I\subseteq(2,3)$ a compact interval, and $\mathcal{N}$ a neighborhood of $(\half,\half)$, such that $\mathcal{N}\subseteq\mathcal{T}_y$ for every $y\in I$. Then there is
\begin{equation*}
\inf\limits_{y\in I}\inf\limits_{(r,\alpha)\in\mathcal{W}\setminus\mathcal{N}}h_y(r,\alpha)>0.
\end{equation*}
\end{corollary}

\begin{proof}
With respect to the function $\alpha_y(r)$ from Lemma~\ref{lem_alphay}, let us introduce the set
\begin{equation*}
\mathdf{R_y}:=\big\{r\in[\mathsmaller{\frac{1}{3}},1)\;\big|\;(r,\alpha_y(r))\notin\mathcal{N}\big\}.
\end{equation*}

\begin{figure}[H]
\centering
\begin{tikzpicture}[scale=1.1]
\draw[->] (1.33,0)--(4.5,0) node[anchor=south] {$r$};
\draw[->] (1.33,0)--(1.33,4.3) node[anchor=west] {$\alpha$};
\fill[black!25] (4,4)--(3.08,4)--(1.33,1.73)--(1.33,0.45);
\draw[ultra thick] (4,4)--(1.33,4)--(1.33,0)--(1.78,0);
\draw[ultra thick] (2.2,0)--(4,0);
\draw[dashed, ultra thick] (4,0)--(4,4);
\draw[ultra thick] (1.33,1.73)--(3.08,4) node[anchor=south] {\tiny{$\frac{2}{y}$}};
\draw[ultra thick] (1.33,0.45)--(4,4);
\draw (3,1) node[anchor=south] {\large{$\cW$}};
\draw (2.8,3.1) node[anchor=north east] {\large{$\mathcal{T}_y$}};
\draw (3.1,3.1) node[anchor=south] {\scriptsize{$\alpha_y(r)$}};
\draw[ultra thick] plot[domain=0.34:0.68] (4*\x-0.02,{(0.067-0.133*\x+\x)*4});
\draw[ultra thick] plot[domain=0.68:0.81] (4*\x-0.02,{(0.4545*(\x-0.81)^2-0.031+\x)*4});
\draw[ultra thick] plot[domain=0.81:1] (4*\x-0.02,{(0.859*(\x-0.81)^2-0.031+\x)*4});
\draw (1.33,0) node[anchor=north] {$\frac{1}{3}$};
\draw (4,0) node[anchor=north] {$1$};
\draw (1.333,0) node[anchor=east] {$0$};
\draw (1.333,4) node[anchor=east] {$1$};
\fill[white] (2,2) circle (0.28cm);
\draw[black, ultra thick] (2,2) circle (0.28cm);
\draw (2,2) node {\small{$\cN$}};
\draw[thick,fill=white] (4,4) circle (0.06cm);
\draw[thick,fill=white] (4,0) circle (0.06cm);
\draw[dashed] (1.78,1.75)--(1.78,0);
\draw[dashed] (2.2,2.23)--(2.2,0);
\draw[ultra thick] (1.78,-0.05)--(1.78,0.05);
\draw[ultra thick] (2.2,-0.05)--(2.2,0.05);
\draw (2,-0.2) node[anchor=north] {\large{$R_y$}};
\draw[->] (1.85,-0.3)--(1.55,-0.03);
\draw[->] (2.15,-0.3)--(3.1,-0.03);
\end{tikzpicture}
\caption{This figure illustrates the neighborhood $\mathcal{N}$ of $(\frac{1}{2},\frac{1}{2})$, as well as the resulting one dimensional set $R_y$ of $r$ values, for which $(r,\alpha_y(r))$ does not intersect with $\mathcal{N}$. The infimum of $h_y(r,\alpha)$ over values $(r,\alpha)\in\mathcal{T}_y\setminus\mathcal{N}$ is now located either on the curve $(r,\alpha_y(r))$, $r\in R_y$ or on the boundary $\partial\mathcal{N}$.}
\end{figure}
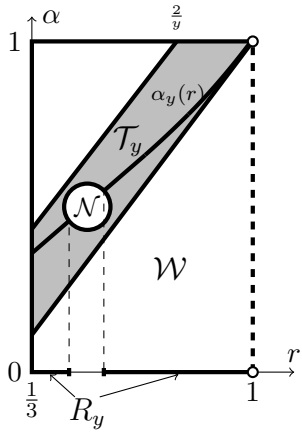

For each $r\in {[\mathsmaller{\frac{1}{3}},1)}$, let us now distinguish two cases \medskip

\begin{itemize}
\item If $(r,\alpha_y(r))\notin\cN$, then $r\in R_y$, and since $\alpha_y(r)$ is the minimizer of the minimization problem \eqref{Eq_hy_Minimization}, we get
\begin{equation}\label{Eq_hy_liminf_2}
\inf\limits_{\alpha\in[0,1]:(r,\alpha)\in\mathcal{T}_y\setminus\mathcal{N}}h_y(r,\alpha)\geq h_y(r,\alpha_y(r))\geq\inf\limits_{\rho\in R_y}h_y(\rho,\alpha_y(\rho)).
\end{equation}

\item If $(r,\alpha_y(r))\in\mathcal{N}$, then the convexity result in Lemma~\ref{lem_hy_convex} shows that $\alpha\mapsto h_y(r,\alpha)$ is minimized for $\alpha$ such that $(r,\alpha)\in\partial\mathcal{N}$ is on the boundary, i.e.,
\begin{equation}\label{Eq_hy_liminf_3}
\inf\limits_{\alpha\in[0,1]:(r,\alpha)\in\mathcal{T}_y\setminus\mathcal{N}}h_y(r,\alpha)=\inf\limits_{\alpha\in[0,1]: (r,\alpha)\in\partial\mathcal{N}}h_y(r,\alpha)\geq\inf\limits_{(\rho,\alpha)\in\partial\mathcal{N}}h_y(\rho,\alpha).
\end{equation}
\end{itemize}

Since the right hand sides of \eqref{Eq_hy_liminf_2} and \eqref{Eq_hy_liminf_3} no longer depends on the chosen value $r\in{[}\frac{1}{3},1)$, we can combine the two cases and get the lower bound of the infimum
\begin{equation}\label{Eq_hy_liminf_4}
\inf\limits_{(r,\alpha)\in\mathcal{T}_y\setminus\mathcal{N}}h_y(r,\alpha)\geq\Big\{\inf\limits_{\rho\in R_y}h_y(\rho,\alpha_y(\rho)),\inf\limits_{(\rho,\alpha)\in\partial\mathcal{N}}h_y(\rho,\alpha)\Big\}.
\end{equation}
For the first infimum on the right hand side of \eqref{Eq_hy_liminf_4}, we note that by \Cref{prop:hy-on-alpha-y}, the only minimizer of the function $h_y(\rho,\alpha_y(\rho))$ is $\rho=\frac{1}{2}$. However, $\frac{1}{2}\notin R_y$, and hence the function $h_y(\rho,\alpha_y(\rho))$ has no local minimum on $R_y$. This means, the minimum is either located at the endpoints $\rho=\frac{1}{3}$ or $\rho\rightarrow 1$, or at the intersection points $(\rho,\alpha_y(\rho))$ with $\partial\mathcal{N}$. Since the latter case is already covered  by the second infimum on the right hand side of Eq.\eqref{Eq_hy_liminf_4}, the inequality \eqref{Eq_hy_liminf_4} becomes
\begin{align*}
\inf\limits_{(r,\alpha)\in\mathcal{T}_y\setminus\mathcal{N}}h_y(r,\alpha)&\geq\min\Big\{h_y(\mathsmaller{\frac{1}{3}},\alpha_y(\mathsmaller{\frac{1}{3}})),\;\lim\limits_{\rho\to 1}h_y(\rho,\alpha_y(\rho)),\inf\limits_{(\rho,\alpha)\in\partial\mathcal{N}}h_y(\rho,\alpha)\bigg\} \\
&\geq\min\Big\{y\big(\mathsmaller{\frac{2\ln(3)}{3}}-\ln(2)\big),\;\ln(2)(1-\mathsmaller{\frac{y}{3}}),\inf\limits_{(\rho,\alpha)\in\partial\mathcal{N}}h_y(\rho,\alpha)\Big\},
\end{align*}
where in the second inequality we used the estimates from Lemma~\ref{lem:hy-near-1}~i)~\&~ii). Finally, using also that $h_y(r,\alpha)=1$, for $(r,\alpha)\in\mathcal{W}\setminus\mathcal{T}_y$ outside the trapezoid, we get
\begin{equation*}
\inf\limits_{(r,\alpha)\in\mathcal{W}\setminus\mathcal{N}}h_y(r,\alpha)\geq\min\Big\{1,\;y\big(\mathsmaller{\frac{2\ln(3)}{3}}-\ln(2)\big),\;\ln(2)(1-\mathsmaller{\frac{y}{3}}),\inf\limits_{(\rho,\alpha)\in\partial\mathcal{N}}h_y(\rho,\alpha)\Big\}.
\end{equation*}
This lower bound for every individual value $y\in(2,3)$ will now be made uniform. From \Cref{thm_hy_Minimum} we know that $\inf_{(\rho,\alpha)\in\partial\mathcal{N}}h_y(r,\alpha)>0$ for every $y\in(2,3)$. Since the mapping $y\mapsto h_y(\rho,\alpha)$ is continuous and $\partial\mathcal{N}$ is compact, also $y\mapsto\inf_{(\rho,\alpha)\in\partial\mathcal{N}}h_y(r,\alpha)$ is continuous. the infimum over the compact interval $I$ then gives
\begin{equation*}
\inf\limits_{y\in I}\inf\limits_{(\rho,\alpha)\in\partial\mathcal{N}}h_y(r,\alpha)=:\delta>0.
\end{equation*}
If we further denote the compact interval $I\subseteq(2,3)$ by $I=[\underline{y},\overline{y}]$, for some $2<\underline{y}\leq\overline{y}<3$, we obtain the uniform lower bound of the infimum
\begin{equation*}
\inf\limits_{y\in I}\inf\limits_{(r,\alpha)\in\mathcal{W}\setminus\mathcal{N}}h_y(r,\alpha)\geq\min\Big\{1,\;\underline{y}\big(\mathsmaller{\frac{2\ln(3)}{3}}-\ln(2)\big),\;\ln(2)\big(1-\mathsmaller{\frac{\overline{y}}{3}}\big),\;\delta\Big\}>0. \qedhere
\end{equation*}
\end{proof}

\subsection{Concluding the proof of \Cref{thm_Application_3XORSAT_appendix}}\label{sec_Proof}

This subsection will finish the proof of \Cref{thm_Application_3XORSAT_appendix}. In \Cref{sec_Subsection1} and \Cref{sec_Subsection2} we verify the assumptions of \Cref{thm_Discrete_Sk} for the summands $\mathcal{S}_{m_k,n_k}$ defined on the smaller domain $\mathcal{W}\cap\Lambda_{m_k,n_k}$ instead of the domain $\widehat{\mathcal{W}}\cap\Lambda_{m_k,n_k}$, which is used in \Cref{thm_Application_3XORSAT_appendix}. In \Cref{sec_Subsection3}, we apply the statement of \Cref{thm_Discrete_Sk}, and also handle the difference between the domains $\mathcal{W}$ and $\widehat{\mathcal{W}}$, to get the explicit value of the limit \eqref{Eq_Application_3XORSAT_appendix}.

\subsubsection{Verifying the \Cref{ass_Laplace_assumption_discrete}}\label{sec_Subsection1}

First we have to check if the function $g_{m_k,n_k}$ and $h_{3m_k/n_k}$ in Eq.\eqref{Eq_gmn} and Eq.\eqref{Eq_hy} satisfy the four parts i)--iv) of Assumption \ref{ass_Laplace_assumption_discrete} with the point $x_0=(\half,\half)$. We will use the value of the limit
\begin{equation*}
\mathdf{c}:=\lim_{k\to\infty}\frac{m_k}{n_k}\in(\mathsmaller{\frac{2}{3}},1),
\end{equation*}
given in the statement of Theorem~\ref{thm_Application_3XORSAT_appendix}. Note that this limit value is in the interval $(\frac{2}{3},1)$ and hence connected to the value $y\in(2,3)$, which appears as an index for the functions $g$ and $h$, in the sense $y=3c$. \medskip

i)\;\;With the function $P$ in Eq.\eqref{Eq_P}, we define for every $y\in(2,3)$ the function
\begin{equation}\label{Eq_g}
\mathdf{g_y(r,\alpha)}\coloneqq\frac{\sqrt{r}\,P(y)}{\pi\alpha(1-\alpha)\sqrt{3r-1}\,P(\frac{ry}{\alpha})P(\frac{(1-r)y}{1-\alpha})}, \qquad (r,\alpha) \in \Int \cT_y.
\end{equation}
Using the explicit representations of $g^S$ and $g^B$ in Eq.\eqref{Eq_gS} and Eq.\eqref{Eq_gB}, we can show that on the interior of $\mathcal{T}_{\frac{3m}{n}}$, the function $g_{m,n}$ in Eq.\eqref{Eq_gmn} coincides with
\begin{equation}\label{Eq_gmn_equals_gy}
g_{m,n}(r,\alpha)=g_{3m/n}(r,\alpha),\qquad(r,\alpha)\in\Int(\cT_{3m/n})\cap\Lambda_{m,n}.
\end{equation}
Since $\lim_{k\to\infty}\frac{3m_k}{n_k}\in(2,3)$, we can choose $N\in\mathbb{N}$ large enough, and a neighborhood $\mathcal{N}_0$ of $(\half,\half)$, such that $\mathcal{N}_0\subseteq\mathcal{T}_{3m_k/n_k}$, for every $k\geq N$, see also the construction in Eq.\eqref{Eq_S2_approximation_1}. Since the mapping $y\mapsto g_y(r,\alpha)$ is continuous for every fixed $(r,\alpha)\in\mathcal{N}_0$, and since the function $g_y$ is continuous for every $y$, there converges
\begin{align*}
\lim\limits_{k\to\infty}\Vert g_{m_k,n_k}-g_{3c}\Vert_{\infty,\mathcal{N}_0\cap\Lambda_{m_k,n_k}}&=\lim\limits_{k\to\infty}\Vert g_{3m_k/n_k}-g_{3c}\Vert_{\infty,\mathcal{N}_0\cap\Lambda_{m_k,n_k}} \\
&\leq\lim\limits_{k\to\infty}\Vert g_{3m_k/n_k}-g_{3c}\Vert_{C(\mathcal{N}_0)}=0.
\end{align*}
Moreover, since also $y\mapsto h_y(r,\alpha)$ is continuous for every $(r,\alpha)\in\mathcal{N}_0$, and $h_y$ is continuous, and since the same is true for all derivatives of $h_y$ up to order two, there also converges
\begin{equation*}
\lim\limits_{k\to\infty}\Vert h_{3m_k/n_k}-h_{3c}\Vert_{C^2(\mathcal{N}_0)}=0.
\end{equation*}

ii)\;\;By \Cref{thm_hy_Minimum}, for every $y\in(2,3)$, the value of $h_y$ at the point $(r,\alpha)=(\half,\half)$ is given by $h_y(\half,\half)=0$. Since this is a local minimum, the gradient $\nabla h_y(\frac{1}{2},\frac{1}{2})=(0,0)$ also vanishes. Moreover, it is straightforward to check that the Hessian matrix in this point is explicitly given by
\begin{equation}\label{Eq_Hy}
\mathdf{\mathcal{H}\{h_y\}(\mathsmaller{\frac{1}{2},\frac{1}{2}})}=\frac{4y^2}{P(y)^2}\begin{pmatrix} 1 & -1 \\ -1 & 1 \end{pmatrix}+\begin{pmatrix} 0 & 0 \\ 0 & 4 \end{pmatrix},
\end{equation}
and hence is clearly positive definite. \medskip

iii)\;\;Since $\lim_{k\to\infty}\frac{3m_k}{n_k}=c\in(\frac{2}{3},1)$, there exists some compact interval $I\subseteq (2,3)$ such that $\frac{3m_k}{n_k} \in I$ for sufficiently large $k$. For every neighborhood $\mathcal{N}$ if $(\frac{1}{2},\frac{1}{2})$ it is now proven in \Cref{cor_hy_liminf}, that
\begin{equation*}
\liminf\limits_{k\to\infty}\inf\limits_{(r,\alpha)\in\cW\setminus\mathcal{N}}h_{3m_k/n_k}(r,\alpha)\geq\inf\limits_{y\in I}\inf\limits_{(r,\alpha)\in\cW\setminus\mathcal{N}}h_y(r,\alpha)>0.
\end{equation*}

iv)\;\;Using  the lower bound of $g^S$ in \eqref{Eq_gS_upper_lower_bound} as well as the upper bound in \eqref{Eq_gS_estimate}, we can upper-bound $g_{m_k,n_k}(r,\alpha)$ defined in \Cref{Eq_gmn} by
\begin{equation*}
g_{m_k,n_k}(r,\alpha)\leq n_k\sqrt{6\pi m_k}\,\frac{g^B(\alpha n_k,n_k)g^B(\frac{3}{2}(1-r)m_k,m_k)}{g^B(3rm_k,3m_k)},\qquad(r,\alpha)\in\mathcal{W}\cap\Lambda_{m_k,n_k}.
\end{equation*}
Using also $\max\{p,1\} \geq 1$ and $\max\{q-p,1\}\geq 1$, yields $g^B(p,q)\leq\frac{\sqrt{q}}{\sqrt{2\pi}}$ for the function $g^B$ in \Cref{Eq_gB}. With this, we can further estimate $g_{m_k,n_k}(r,\alpha)$ by
\begin{align*}
g_{m_k,n_k}(r,\alpha)&\leq\frac{\sqrt{3}\,n_k^{\frac{3}{2}}m_k}{\sqrt{2\pi}\,g^B(3rm_k,3m_k)} \\
&=3(n_km_k)^{\frac{3}{2}}\sqrt{r(1-r)}\leq\frac{3}{2}C^{\frac{3}{2}}n_k^3,\qquad(r,\alpha)\in\mathcal{W}\cap\Lambda_{m_k,n_k},
\end{align*}
where in the second line we plugged in the definition \eqref{Eq_gB} of the function $g_B$, and in the third inequality, we used $\sqrt{r(1-r)}\leq\frac{1}{2}$ as well as the constant $C:=\sup_{k\in\mathbb{N}}\frac{m_k}{n_k}<\infty$.

\subsubsection{Verifying the assumptions of \Cref{thm_Discrete_Sk}}\label{sec_Subsection2}

Next, we will verify that the summands $\mathcal{S}_{m_k,n_k}$ in \eqref{Eq_Smn_appendix} satisfy the assumptions i)~\&~ii) of Theorem \ref{thm_Discrete_Sk}. \medskip

i)\;\;It is already shown in Eq.\eqref{Eq_Smn_estimate}, that
\begin{equation*}
\mathcal{S}_{m_k,n_k}(r,\alpha)\leq\frac{C_S^2C_B^2}{c_Sc_B}g_{m_k,n_k}(r,\alpha)e^{-n_kh_{3m_k/n_k}(r,\alpha)},\qquad(r,\alpha)\in\mathcal{W}\cap\Lambda_{m_k,n_k}.
\end{equation*}

ii)\;\;As in \Cref{sec_Subsection1}~iii), there exists $2<\underline{y}\leq\overline{y}<3$ and $N\in\mathbb{N}$, such that $\frac{3m_k}{n_k}\in[\underline{y},\overline{y}]$, for every $k\geq N$. Moreover, we can choose a neighborhood $\mathcal{N}_0$ of $(\half,\half)$ small enough, such that $\mathcal{N}_0\subseteq\mathcal{T}_{\underline{y}}$, and there is
\begin{equation*}
\alpha(1-\alpha)\geq\frac{1}{8},\qquad 3r-1\geq\frac{1}{4},\qquad\frac{2}{3\overline{c}}\leq\frac{r}{\alpha},\qquad\text{and}\qquad\frac{2}{3\overline{c}}\leq\frac{1-r}{1-\alpha},\qquad(r,\alpha)\in\mathcal{N}_0.
\end{equation*}
With the representations \eqref{Eq_g} and \eqref{Eq_gmn_equals_gy}, we can now use the monotonicity of the function $P$, to estimate $g_{m_k,n_k}(r,\alpha)$ for every $(r,\alpha)\in\mathcal{N}_0\cap\Lambda_{m_k,n_k}$ by
\begin{equation*}
g_{m_k,n_k}(r,\alpha)=\frac{\sqrt{r}\,P(\frac{3m_k}{n_k})}{\pi\alpha(1-\alpha)\sqrt{3r-1}\,P(\frac{3rm_k}{\alpha n_k})P(\frac{3(1-r)m_k}{(1-\alpha)n_k})}\leq\frac{16P(\overline{y})}{\pi P(\frac{2\underline{y}}{\overline{y}})^2}\eqqcolon\mathdf{C_{\mathcal{N}_0}}.
\end{equation*}
With this uniform upper bound of $g_{m_k,n_k}$, we
can estimate for every $(r,\alpha)\in\mathcal{N}_0\cap\Lambda_{n_k,m_k}$
\begin{equation*}
\Big|\mathcal{S}_{m_k,n_k}(r,\alpha)e^{n_kh_{3m_k/n_k}(r,\alpha)}-g_{m_k,n_k}(r,\alpha)\Big|\leq C_{\mathcal{N}_0}\bigg|\frac{\mathcal{S}_{m_k,n_k}(r,\alpha)}{g_{m_k,n_k}(r,\alpha)e^{-n_kh_{3m_k/n_k}(r,\alpha)}}-1\bigg|.
\end{equation*}
Since we already know that the asymptotics \eqref{Eq_Smn_asymptotics} is uniform on $\mathcal{N}_0\cap\Lambda_{m_k,n_k}$, if we make $\mathcal{N}_0$ smaller if necessary, we conclude the uniform convergence
\begin{align*}
\lim\limits_{k\to\infty}\big\Vert\mathcal{S}_{m_k,n_k}&e^{n_kh_{3m_k/n_k}}-g_{m_k,n_k}\big\Vert_{\infty,\mathcal{N}_0\cap\Lambda_{m_k,n_k}} \\
&\leq C_{\mathcal{N}_0}\lim\limits_{k\to\infty}\bigg\Vert\frac{\mathcal{S}_{m_k,n_k}}{g_{m_k,n_k}e^{-n_kh_{3m_k/n_k}}}-1\bigg\Vert_{\infty,\mathcal{N}_0\cap\Lambda_{m_k,n_k}}=0.
\end{align*}

\subsubsection{Applying \Cref{thm_Discrete_Sk}}\label{sec_Subsection3}

Note that above we verified the hypotheses of
\Cref{thm_Discrete_Sk} for 
$\cW=[\frac{1}{3},1)\times[0,1]$, 
but \Cref{thm_Application_3XORSAT_appendix} makes a claim about a summation over $\hW=[\frac{1}{3},1]\times[0,1]$. However, the difference is easily handled. Indeed, there is $\mathcal{S}_{m_k,n_k}(r,\alpha) = 0$ for $r=1$, except at $(r,\alpha)=(1,1)$ where it can be easily checked that $\mathcal{S}_{m_k,n_k}(1,1) = n_k 2^{m_k-n_k}$. Hence the two sums over $\widehat{\mathcal{W}}\cap\Lambda_{m_k,n_k}$ and $\mathcal{W}\cap\Lambda_{m_k,n_k}$ only differ by
\begin{align*}
\frac{1}{n_k}\sum\limits_{(r,\alpha)\in\hW\cap\Lambda_{m_k,n_k}}\mathcal{S}_{m_k,n_k}(r,\alpha)
= 2^{m_k-n_k} + \frac{1}{n_k}\sum\limits_{(r,\alpha)\in\cW\cap\Lambda_{m_k,n_k}}\mathcal{S}_{m_k,n_k}(r,\alpha).
\end{align*}

Because $\lim_{k\to\infty}\frac{m_k}{n_k}\in(\frac{2}{3},1)$, we get $\lim_{k\to\infty}2^{m_k-n_k}=0$, and so
\begin{align*}
\lim\limits_{k\to\infty}\frac{1}{n_k}\sum\limits_{(r,\alpha)\in\hW\cap\Lambda_{m_k,n_k}}\mathcal{S}_{m_k,n_k}(r,\alpha)
=
\lim\limits_{k\to\infty}\frac{1}{n_k}\sum\limits_{(r,\alpha)\in\cW\cap\Lambda_{m_k,n_k}}\mathcal{S}_{m_k,n_k}(r,\alpha)
\end{align*}

In the Sections~\ref{sec_Subsection1}~\&~\ref{sec_Subsection2} we have verified that all the assumptions of Theorem \ref{thm_Discrete_Sk} are satisfied for $\cW$, and we conclude the limit

\begin{equation}\label{Eq_Final_limit_appendix}
\lim\limits_{k\to\infty}\frac{1}{n_k}\sum\limits_{(r,\alpha)\in\hW\cap\Lambda_{m_k,n_k}}\mathcal{S}_{m_k,n_k}(r,\alpha)=\frac{2\pi g_{3c}(\frac{1}{2},\frac{1}{2})}{|\det(A)|\sqrt{\det(\mathcal{H}\{h_{3c}\}(\frac{1}{2},\frac{1}{2}))}}.
\end{equation}
Here, the matrix $A$ comes from the decomposition \eqref{Eq_Lambdamn} of the lattice $\Lambda_{m_k,n_k}$ in the sense \eqref{Eq_Lattice}, i.e.,
\begin{equation*}
\Lambda_{m_k,n_k}=\frac{1}{n_k}\underbrace{\begin{pmatrix} \frac{2n_k}{3m_k} & 0 \\ 0 & 1 \end{pmatrix}}_{=:A_k}+\begin{pmatrix} 1 \\ 0 \end{pmatrix},
\end{equation*}
and the resulting convergence
\begin{equation}\label{Eq_A_matrix}
\lim\limits_{k\to\infty}A_k=\lim\limits_{k\to\infty}\begin{pmatrix} \frac{2n_k}{3m_k} & 0 \\ 0 & 1 \end{pmatrix}=\begin{pmatrix} \frac{2}{3c} & 0 \\ 0 & 1 \end{pmatrix}=:A.
\end{equation}
The exact values of the terms on the right hand side of Eq.\eqref{Eq_Final_limit_appendix} can now be calculated from their explicit representations in \eqref{Eq_g}, \eqref{Eq_A_matrix} and \eqref{Eq_Hy}
\begin{align*}
g_{3c}(\mathsmaller{\frac{1}{2},\frac{1}{2}})=\frac{4}{\pi\,P(3c)},\qquad\det(A)=\frac{2}{3c},\qquad\det\big(\mathcal{H}\{h_{3c}\}(\mathsmaller{\frac{1}{2},\frac{1}{2}})\big)=\frac{16(3c)^2}{P(3c)^2}.
\end{align*}
With these values we then conclude the final limit
\begin{equation*}
\lim\limits_{k\to\infty}\frac{1}{n_k}\sum\limits_{(r,\alpha)\in\hW\cap\Lambda_{m_k,n_k}}\mathcal{S}_{m_k,n_k}(r,\alpha)=1.
\end{equation*}
This concludes the proof of Theorem~\ref{thm_Application_3XORSAT_appendix}. \hfill\qed

\vspace{1cm}

\clearpage

\label{table-of-contents}
\hypertarget{table-of-contents}{}
\setcounter{tocdepth}{3}
\tableofcontents

\clearpage
\printindex

\end{document}